\documentclass[aos,MSNbibl,seceqn,dvips]{arximspdf}
\usepackage{graphicx}
\doi{10.1214/14-AOS1267} 
\volume{43}
\issue{1}
\pubyear{2015}
\firstpage{102}
\lastpage{138}
\docsubty{FLA}

\makeatletter
\newcommand{\T}{\intercal}
\renewcommand{\mid}{|}
\newcommand{\rrvert}{\vert}
\newcommand{\rrVert}{\Vert}
\newcommand{\llvert}{\vert}
\newcommand{\llVert}{\Vert}
\newproclaim{Definition}{Definition}[section]
\newtheorem{Theorem}{Theorem}[section]
\newtheorem{Lemma}{Lemma}[section]
\newproclaim{Remark}{Remark}[section]
\newtheorem{Corollary}{Corollary}[section]
\newtheorem{Proposition}{Proposition}[section]
\makeatother

\begin{document}
\begin{frontmatter}

\title{ROP: Matrix recovery via rank-one projections\thanksref{T1}}
\runtitle{Rank-one projections}

\begin{aug}
\author[A]{\fnms{T. Tony}~\snm{Cai}\corref{}\ead[label=e1]{tcai@wharton.upenn.edu}\ead[label=u1,url]{http://www-stat.wharton.upenn.edu/\textasciitilde tcai/}}
\and
\author[A]{\fnms{Anru}~\snm{Zhang}\ead[label=e2]{anrzhang@wharton.upenn.edu}\ead[label=u2,url]{http://www-stat.wharton.upenn.edu/\textasciitilde anrzhang/}}
\runauthor{T.~T. Cai and A. Zhang}
\affiliation{University of Pennsylvania}
\address[A]{Department of Statistics\\
The Wharton School\\
University of Pennsylvania\\
Philadelphia, Pennsylvania 19104\\
USA\\
\printead{e1}\\
\phantom{E-mail: }\printead*{e2}\\
\printead{u1}\\
\phantom{URL: }\printead*{u2}}
\end{aug}
\thankstext{T1}{Supported in part by NSF FRG Grant DMS-08-54973,
NSF Grant DMS-12-08982 and NIH Grant R01 CA127334-05.}

\received{\smonth{3} \syear{2014}}
\revised{\smonth{8} \syear{2014}}

%
\begin{abstract}
Estimation of low-rank matrices is of significant interest in a range
of contemporary applications.
In this paper, we introduce a rank-one projection model for low-rank
matrix recovery and propose
a constrained nuclear norm minimization method for stable recovery of
low-rank matrices in the noisy
case. The procedure is adaptive to the rank and robust against small
perturbations. Both upper and lower bounds for the estimation accuracy
under the Frobenius norm loss are obtained. The proposed estimator is
shown to be rate-optimal under certain conditions. The estimator is
easy to implement via convex programming and performs well numerically.

The techniques and main results developed in the paper also have
implications to other related statistical
problems. An application to estimation of spiked covariance matrices
from one-dimensional random projections is considered. The results
demonstrate that it is still possible to accurately estimate the
covariance matrix of a high-dimensional distribution based only on
one-dimensional projections.
\end{abstract}

%
\end{frontmatter}

\section{Introduction}\label{Introductionsec}

Accurate recovery of low-rank matrices has a wide range of
applications, including quantum state tomography \cite{Alquier,Gross},
face recognition \mbox{\cite{Basri,Candes-Li}}, recommender systems \cite
{Koren} and linear system identification and control \cite
{RechtMatrix}. For example, a key step in reconstructing the quantum
states in low-rank quantum tomography is the estimation of a low-rank
matrix based on Pauli measurements \cite{Gross,Wang}. And phase
retrieval, a problem which arises in a range of signal and image
processing applications including X-ray crystallography, astronomical
imaging and diffraction imaging, can be reformulated as a low-rank
matrix recovery problem \cite{CandesRIPfail,Candes-Li}. See Recht et~al. \cite{RechtMatrix} and Cand\`es and Plan \cite{CandesOracle}
for further references and discussions.

Motivated by these applications, low-rank matrix estimation based on a
small number of measurements has drawn much recent attention in several
fields, including statistics, electrical engineering, applied
mathematics and computer science. For example, Cand\`es and Recht \cite
{CandesRecht}, Cand\`es and Tao \cite{CandesTao} and Recht \cite
{RechtImproved} considered the exact recovery of a low-rank matrix
based on a subset of uniformly sampled entries. Negahban and Wainwright
\cite{Negahban} investigated matrix completion under a row/column
weighted random sampling scheme. Recht et~al. \cite{RechtMatrix},
Cand\`es and Plan \cite{CandesOracle} and Cai and Zhang \cite
{CaiZhang,CaiZhang2,CaiZhang3} studied matrix recovery based on a
small number of linear measurements in the framework of Restricted Isometry Property (RIP),
and Koltchinskii et~al. \cite{KLT} proposed the penalized nuclear norm
minimization method and derived a general sharp oracle inequality under
the condition of restrict isometry in expectation.

The basic model for low-rank matrix recovery can be written as
%
\begin{equation}
\label{eqmodel} y = \mathcal{X}(A)+z,
\end{equation}
where $\mathcal{X}\dvtx  \mathbb{R}^{p_1\times p_2} \to\mathbb{R}^n$ is
a linear map, $A\in\mathbb{R}^{p_1\times p_2}$ is an unknown low-rank
matrix and $z$ is a noise vector. The goal is to recover the low-rank
matrix $A$ based on the measurements $(\mathcal{X}, y)$. The linear
map $\mathcal{X}$ can be equivalently specified by $n$ $p_1\times p_2$
measurement matrices $X_1,\ldots, X_n$ with
%
\begin{equation}
\label{eqM} \mathcal{X}(A) = \bigl(\langle X_1, A\rangle, \langle
X_2, A\rangle, \ldots, \langle X_n, A\rangle
\bigr)^\intercal,
\end{equation}
where the inner product of two matrices of the same dimensions is
defined as $\langle X, Y\rangle= \sum_{i,j} X_{ij}Y_{ij}$. Since
$\langle X, Y\rangle=\operatorname{trace}(X^\intercal Y)$, (\ref{eqmodel}) is
also known as trace regression.

A common approach to low-rank matrix recovery is the constrained
nuclear norm minimization method which estimates $A$ by
%
\begin{equation}
\label{eqnuclearnormminimization} \hat A = \mathop{\arg\min}_{M} \bigl\{\llVert M \rrVert
_\ast\dvtx  y-\mathcal{X}(M) \in \mathcal{Z} \bigr\}.
\end{equation}
Here, $\llVert  X \rrVert  _\ast$ is the nuclear norm of the matrix $X$ which is
defined to be the sum of its singular values, and $\mathcal{Z}$ is a
bounded set determined by the noise structure. For example, $\mathcal
{Z} = \{0\}$ in the noiseless case and $\mathcal{Z}$ is the feasible
set of the error vector $z$ in the case of bounded noise. This
constrained nuclear norm minimization method has been well studied.
See, for example, \cite
{RechtMatrix,CandesOracle,Oymak,CaiZhang,CaiZhang2,CaiZhang3}.

Two random design models for low-rank matrix recovery have been
particularly well studied in the literature. One is the so-called
``Gaussian ensemble'' \cite{RechtMatrix,CandesOracle}, where the
measurement matrices $X_1, \ldots, X_n$ are random matrices with
i.i.d. Gaussian entries. By exploiting the low-dimensional structure,
the number of linear measurements can be far smaller than the number of
entries in the matrix to ensure stable recovery. It has been shown that
a matrix $A$ of rank $r$ can be stably recovered by nuclear norm
minimization with high probability, provided that $n \gtrsim
r(p_1+p_2)$ \cite{CandesOracle}. One major disadvantage of the
Gaussian ensemble design is that it requires $O(np_1p_2)$ bytes of
storage space for $\mathcal{X}$, which can be excessively large for
the recovery of large matrices. For example, at least $45$ TB of space
is need to store the measurement matrices $M_i$ in order to ensure
accurate reconstruction of 10,000${}\times{}$10,000 matrices of rank 10.
(See more discussion in Section~\ref{Simulationssec}.) Another
popular design is the ``matrix completion'' model \cite
{CandesRecht,CandesTao,RechtImproved}, under which the individual
entries of the matrix $A$ are observed at randomly selected positions.
In terms of the measurement matrices $X_i$ in (\ref{eqM}), this can
be interpreted as
%
\begin{equation}
\label{MC} \mathcal{X}(A) = \bigl(\bigl\langle e_{i_1}e_{j_1}^\intercal,
A\bigr\rangle, \bigl\langle e_{i_2}e_{j_2}^\intercal, A
\bigr\rangle, \ldots, \bigl\langle e_{i_n}e_{j_n}^\intercal,
A\bigr\rangle\bigr)^\intercal,
\end{equation}
where $e_i = (0, \ldots, 0, \overbrace{1}^{i {\mathrm{th}}}, 0,
\ldots, 0)$ is the $i$th standard basis vector, and $i_1,\ldots,\break  i_n$
and $j_1, \ldots, j_n$ are randomly and uniformly drawn with
replacement from $\{1, \ldots, p_1\}$ and $\{1, \ldots, p_2\}$,
respectively. However, as pointed out in \cite
{CandesRecht,RechtImproved}, additional structural assumptions, which
are not intuitive and difficult to check, on the unknown matrix $A$ are
needed in order to ensure stable recovery under the matrix completion
model. For example, it is impossible to recover spiked matrices under
the matrix completion model. This can be easily seen from a simple
example where the matrix $A$ has only one nonzero row. In this case,
although the matrix is only of rank one, it is not recoverable under
the matrix completion model unless all the elements on the nonzero row
are observed.

In this paper, we introduce a ``\emph{Rank-One Projection}'' (\emph{ROP})
model for low-rank matrix recovery and propose a constrained nuclear
norm minimization method for this model. Under the ROP model, we observe
%
\begin{equation}
\label{eqROP-model} y_i =\bigl(\beta^{(i)}\bigr)^\intercal A
\gamma^{(i)} +z_i, \qquad i =1, \ldots, n,
\end{equation}
where $\beta^{(i)}$ and $\gamma^{(i)}$ are random vectors with
entries independently drawn from some distribution $\mathcal{P}$, and
$z_i$ are random errors. In terms of the linear map $\mathcal{X}\dvtx
\mathbb{R}^{p_1 \times p_2} \to\mathbb{R}^n$ in (\ref{eqmodel}),
it can be defined as
%
\begin{equation}
\label{ROP2} \bigl[\mathcal{X}(A)\bigr]_i =\bigl(
\beta^{(i)}\bigr)^{\intercal} A \gamma^{(i)}, \qquad i=1,
\ldots, n.
\end{equation}
Since the measurement matrices $X_i = \beta^{(i)} (\gamma
^{(i)})^\intercal$ are of rank-one, we call the model (\ref
{eqROP-model}) a ``\emph{Rank-One Projection}'' (\emph{ROP}) model.
It is easy to see that the storage for the measurement vectors in the
ROP model (\ref{eqROP-model}) is $O(n(p_1+p_2))$ bytes which is
significantly smaller than $O(np_1p_2)$ bytes required for the Gaussian
ensemble.

We first establish a sufficient identifiability condition in
Section~\ref{Gaussiansec} by considering the problem of exact
recovery of low-rank matrices in the noiseless case. It is shown that,
with high probability, ROP with $n \gtrsim r(p_1+p_2)$ random
projections is sufficient to ensure exact recovery of all rank-$r$
matrices through the constrained nuclear norm minimization. The
required number of measurements $O( r(p_1+p_2))$ is rate optimal for
any linear measurement model since a rank-$r$ matrix $A\in\mathbb
{R}^{p_1+p_2}$ has the degree of freedom $r(p_1 + p_2 - r)$.
The Gaussian noise case is of particular interest in statistics. We
propose a new constrained nuclear norm minimization estimator and
investigate its theoretical and numerical properties in the Gaussian
noise case. Both upper and lower bounds for the estimation accuracy
under the Frobenius norm loss are obtained. The estimator is shown to
be rate-optimal when the number of rank-one projections satisfies
either $n \gtrsim(p_1+p_2)\log(p_1 + p_2)$ or $n \sim r(p_1+p_2)$.
The lower bound also shows that if the number of measurements $n <
r\max(p_1, p_2)$, then no estimator can recover rank-$r$ matrices
consistently. The general case where the matrix $A$ is only
approximately low-rank is also considered. The results show that the
proposed estimator is adaptive to the rank $r$ and robust against small
perturbations. Extensions to the sub-Gaussian design and sub-Gaussian
noise distribution are also considered.

The ROP model can be further simplified by taking $\beta^{(i)}=\gamma
^{(i)}$ if the low-rank matrix $A$ is known to be symmetric. This is
the case in many applications, including low-dimensional Euclidean
embedding \cite{Trosset,RechtMatrix}, phase retrieval \cite
{CandesRIPfail,Candes-Li} and covariance matrix estimation \cite
{Chen,CaiMaWu1,CaiMaWu2}. In such a setting, the ROP design can be
simplified to Symmetric Rank-One Projections (SROP)
\[
\bigl[\mathcal{X}(A)\bigr]_i = \bigl(\beta^{(i)}
\bigr)^\intercal A \beta^{(i)}.
\]
We will show that the results for the general ROP model continue to
hold for the SROP model when $A$ is known to be symmetric. Recovery of
symmetric positive definite matrices in the noiseless and $\ell
_1$-bounded noise settings has also been considered in a recent paper
by Chen et~al. \cite{Chen} which was posted on arXiv at the time of
the writing of the present paper. Their results and techniques for
symmetric positive definite matrices are not applicable to the recovery
of general low-rank matrices. See Section~\ref{Discussionssec} for
more discussions.

The techniques and main results developed in the paper also have
implications to other related statistical problems. In particular, the
results imply that it is possible to accurately estimate a spiked
covariance matrix based only on one-dimensional projections. Spiked
covariance matrix model has been well studied in the context of
Principal Component Analysis (PCA) based on i.i.d. data where one
observes $p$-dimensional vectors $X^{(1)}, \ldots, X^{(n)} \stackrel
{\mathrm{i.i.d.}}{\sim} N(0,\Sigma)$ with $\Sigma= I_p + \Sigma_0$ and $\Sigma
_0$ being low-rank \cite{Johnstone,Birnbaum,CaiMaWu1,CaiMaWu2}. This
covariance structure and its variations have been used in many
applications including signal processing, financial econometrics,
chemometrics and population genetics. See, for example, \cite
{Fan,Nadler,Patterson,Price,Wax}. Suppose that the random vectors
$X^{(1)}, \ldots, X^{(n)}$ are not\vspace*{1pt} directly observable. Instead, we
observe only one-dimensional random projections of $X^{(i)}$,
\[
\xi_{i} = \bigl\langle\beta^{(i)}, X^{(i)} \bigr
\rangle, \qquad i=1, \ldots, n,
\]
where $\beta^{(i)}\stackrel{\mathrm{i.i.d.}}{\sim} N(0, I_p)$. It is somewhat
surprising that it is still possible to accurately estimate the spiked
covariance matrix $\Sigma$ based only on the one-dimensional
projections $\{\xi_i\dvtx  i=1, \ldots, n\}$. This covariance matrix
recovery problem is also related to the recent literature on covariance
sketching \cite{Nowak1,Nowak2}, which aims to recover a symmetric
matrix $A$ (or a general rectangular matrix $B$) from low-dimensional
projections of the form $X^\intercal A X$ (or $X^\intercal B Y$). See
Section~\ref{Applicationssec} for further discussions.

The proposed methods can be efficiently implemented via convex
programming. A simulation study is carried out to investigate the
numerical performance of the proposed nuclear norm minimization
estimators. The numerical results indicate that ROP with $n \ge5 r\max
(p_1, p_2)$ random projections is sufficient to ensure the exact
recovery of rank-$r$ matrices through constrained nuclear norm
minimization and show that the procedure is robust against small
perturbations, which confirm the theoretical results developed in the
paper. The proposed estimator outperforms two other alternative
procedures numerically in the noisy case. In addition, the proposed
method is illustrated through an image compression example.

The rest of the paper is organized as follows. In Section~\ref
{Gaussiansec}, after introducing basic notation and definitions, we
consider exact recovery of low-rank matrices in the noiseless case and
establish a sufficient identifiability condition. A constrained nuclear
norm minimization estimator is introduced for the Gaussian noise case.
Both upper and lower bounds are obtained for estimation under the
Frobenius norm loss. Section~\ref{Extensionssec} considers extensions
to sub-Gaussian design and sub-Gaussian noise distributions. An
application to estimation of spiked covariance matrices based on
one-dimensional projections is discussed in detail in Section~\ref
{Applicationssec}. Section~\ref{Simulationssec} investigates the
numerical performance of the proposed procedure through a simulation
study and an image compression example.
A brief discussion is given in Section~\ref{Discussionssec}. The main
results are proved in Section~\ref{Proofssec} and the proofs of some
technical lemmas are given in the
supplementary material \cite{Supplement}.

\section{Matrix recovery under Gaussian noise}\label{Gaussiansec}

In this section, we first establish an identifiability condition for
the ROP model by considering exact recovery in the noiseless case, and
then focus on low-rank matrix recovery in the Gaussian noise case.

We begin with the basic notation and definitions. For a vector $\beta
\in\mathbb{R}^n$, we use $\llVert  \beta \rrVert  _q = \sqrt[q]{\sum_{i=1}^n\llvert  \beta_i \rrvert  ^q}$ to define its vector $q$-norm. For a matrix
$X\in\mathbb{R}^{p_1\times p_2}$, the Frobenius norm is $\llVert  X \rrVert  _F =
\sqrt{\sum_{i=1}^{p_1}\sum_{j=1}^{p_2}X_{ij}^2}$ and the spectral
norm $\llVert  \cdot \rrVert  $ is $\llVert  X \rrVert   = \sup_{\llVert  \beta \rrVert  _2 \leq1}\llVert  X\beta \rrVert
_2$. For a linear map $\mathcal{X}=(X_1, \ldots, X_n)$ from $\mathbb
{R}^{p_1\times p_2}$ to $\mathbb{R}^n$ given by~(\ref{eqM}), its
dual operator $\mathcal{X}^\ast\dvtx \mathbb{R}^n \to\mathbb
{R}^{p_1\times p_2}$ is defined as $\mathcal{X}^\ast(z) = \sum_{i=1}^n z_i X_i$. For a matrix $X\in\mathbb{R}^{p_1\times p_2}$,
let $X=\sum_{i} a_i u_i v_i^\intercal$ be the singular value
decomposition of $X$ with the singular values $a_1\ge a_2 \ge\cdots
\ge0$. We define $X_{\max(r)}=\sum_{i=1}^ra_iu_iv_i^\intercal$ and
$X_{-\max(r)}=X-X_{\max(r)}=\sum_{i\ge r+1} a_iu_iv_i^\intercal$.
For any two sequences $\{a_n\}$ and $\{b_n\}$ of positive numbers,
denote by $a_n\gtrsim b_n$ when $a_n\geq C b_n$ for some uniform
constant $C$ and denote by $a_n\sim b_n$ if $a_n \gtrsim b_n$ and
$b_n\gtrsim a_n$.

We use the phrase ``rank-$r$ matrices'' to refer to matrices of rank at
most $r$ and denote by $\mathbb{S}^{p}$ the set of all $p\times p$
symmetric matrices.
A linear map $\mathcal{X}\dvtx \mathbb{R}^{p_1\times p_2}\to\mathbb
{R}^n$ is called ROP from distribution $\mathcal{P}$ if $\mathcal{X}$
is defined as in (\ref{ROP2}) with all the entries of $\beta^{(i)}$
and $\gamma^{(i)}$ independently drawn from the distribution $\mathcal{P}$.

\subsection{RUB, identifiability, and exact recovery in the noiseless case}\label{Noiselesssec}

An important step toward understanding the constrained nuclear norm
minimization is the study of exact recovery of low-rank matrices in the
noiseless case which also leads to a sufficient identifiability
condition. A widely used framework in the low-rank matrix recovery
literature is the Restricted Isometry Property (RIP) in the matrix
setting. See \cite
{RechtMatrix,CandesOracle,Rohde,CaiZhang,CaiZhang2,CaiZhang3}.
However, the RIP framework is not well suited for the ROP model and
would lead to suboptimal results. See Section~\ref
{otherconditionssec} for more discussions on the RIP and other
conditions used in the literature. See also \cite{CandesRIPfail}. In
this section, we introduce a Restricted Uniform Boundedness (RUB)
condition which will be shown to guarantee the exact recovery of
low-rank matrices in the noiseless case and stable recovery in the
noisy case through the constrained nuclear norm minimization. It will
also be shown that the RUB condition are satisfied by a range of random
linear maps with high probability.

%
\begin{Definition}[(Restricted Uniform Boundedness)]\label{dfl1bound}
For a linear map $\mathcal{X}\dvtx  \mathbb{R}^{p_1\times p_2} \to\mathbb
{R}^n$, if there exist uniform constants $C_1$ and $ C_2$ such that for
all nonzero rank-$r$ matrices $A\in\mathbb{R}^{p_1\times p_2}$
\[
C_1 \leq{\llVert  \mathcal{X}(A)\rrVert  _1/n\over\llVert  A \rrVert  _F} \leq C_2 ,
\]
where $\llVert  \cdot \rrVert  _1$ means the vector $\ell_1$ norm, then we say that
$\mathcal{X}$ satisfies the Restricted Uniform Boundedness (RUB)
condition of order $r$ and constants $C_1$ and $ C_2$.
\end{Definition}

In the noiseless case, we observe $y = \mathcal{X}(A)$ and estimate
the matrix $A$ through the constrained nuclear norm minimization
%
\begin{equation}
\label{eqnuclearnormminimizationnoiseless} A_\ast= \mathop{\arg\min}_{M} \bigl\{\llVert
M \rrVert _\ast\dvtx  \mathcal{X}(M)= y \bigr\}.
\end{equation}
The following theorem shows that the RUB condition guarantees the exact
recovery of all rank-$r$ matrices.

%
\begin{Theorem}\label{thl1boundrecovery}
Let $k\geq2$ be an integer. Suppose $\mathcal{X}$ satisfies RUB of
order $kr$ with $C_2/C_1< \sqrt{k}$, then the nuclear norm
minimization method recovers all \mbox{rank-$r$} matrices. That is, for all
rank-$r$ matrices $A$ and $y= \mathcal{X}(A)$, we have $A_\ast= A$,
where $A_\ast$ is given by (\ref{eqnuclearnormminimizationnoiseless}).
\end{Theorem}

Theorem \ref{thl1boundrecovery} shows that RUB of order $kr$ with
$C_2/C_1< \sqrt{k}$ is a sufficient identifiability condition for the
low-rank matrix recovery model (\ref{eqmodel}) in the noisy case. The
following result shows that the RUB condition is satisfied with high
probability under the ROP model with a sufficient number of measurements.

\begin{Theorem}\label{thgaussianl1bound}
Suppose\vspace*{1pt} $\mathcal{X}\dvtx \mathbb{R}^{p_1\times p_2} \to\mathbb{R}^n$ is
ROP from the standard normal distribution. For integer $k \geq2$,
positive numbers $C_1<\frac{1}{3}$ and $C_2>1$, there exist constants
$C$ and $\delta$, not depending on $p_1, p_2$ and $r$, such that if
%
\begin{equation}
\label{nlowerbd} n \geq Cr(p_1+p_2),
\end{equation}
then with probability at least $1-e^{-n\delta}$, $\mathcal{X}$
satisfies RUB of order $kr$ and constants $C_1$ and $C_2$.
\end{Theorem}

%
\begin{Remark}
The condition $n \geq O(r(p_1 + p_2))$ on the
number of measurements is indeed necessary for $\mathcal{X}$ to
satisfy nontrivial RUB with $C_1>0$. Note that the degree of freedom of
all rank-$r$ matrices of $\mathbb{R}^{p_1\times p_2}$ is $r(p_1+p_2-r)
\geq\frac{1}{2}r(p_1+p_2)$. If $n < \frac{1}{2}r(p_1 + p_2)$, there
must exist a nonzero rank-$r$ matrix $A\in\mathbb{R}^{p_1\times p_2}$
such that $\mathcal{X}(A)=0$, which leads to the failure of any
nontrivial RUB for $\mathcal{X}$.
\end{Remark}

As a direct consequence of Theorems \ref{thl1boundrecovery} and \ref
{thgaussianl1bound}, ROP with the number of measurements $n \geq
Cr(p_1+p_2)$ guarantees the exact recovery of all rank-$r$ matrices
with high probability.

%
\begin{Corollary}\label{rmGaussiannoiseless}
Suppose $\mathcal{X}\dvtx \mathbb{R}^{p_1\times p_2}\to\mathbb{R}^n$ is
ROP from the standard normal distribution. There exist uniform
constants $C$ and $\delta$ such that, whenever $n \geq Cr(p_1+p_2)$,
the nuclear norm minimization estimator $A_*$ given in (\ref
{eqnuclearnormminimizationnoiseless}) recovers all rank-$r$
matrices $A\in\mathbb{R}^{p_1\times p_2}$ exactly with probability at
least $1-e^{-n\delta}$.
\end{Corollary}

Note that the required number of measurements $O(r(p_1+p_2))$ above is
rate optimal, since the degree of freedom for a matrix $A\in\mathbb
{R}^{p_1+p_2}$ of rank $r$ is $r(p_1 + p_2 - r)$, and thus at least
$r(p_1 + p_2 - r)$ measurements are needed in order to recover $A$
exactly using any method.

\subsection{RUB, RIP and other conditions}\label{otherconditionssec}

We have shown that RUB implies exact recovery in the noiseless and
proved that the random rank-one projections satisfy RUB with high
probability whenever the number of measurements $n \geq Cr(p_1+p_2)$.
As mentioned earlier, other conditions, including the Restricted
Isometry Property (RIP), RIP in expectation and Spherical Section
Property (SSP), have been introduced for low-rank matrix recovery based
on linear measurements.
Among them, RIP is perhaps the most widely used. 
A linear map $\mathcal{X}\dvtx  \mathbb{R}^{p_1\times p_2} \to\mathbb
{R}^{n}$ is said to satisfy RIP of order $r$ with positive constants
$C_1$ and $C_2$ if
\[
C_1 \leq\frac{\llVert  \mathcal{X}(A)\rrVert  _2/\sqrt{n}}{\llVert  A \rrVert  _F} \leq C_2 %
\]
for all rank-$r$ matrices $A$. Many results have been given for
low-rank matrices under the RIP framework. For example, Recht et~al.
\cite{RechtMatrix} showed that Gaussian ensembles satisfy RIP with
high probability under certain conditions on the dimensions. Cand\`es
and Plan \cite{CandesOracle} provided a lower bound and oracle
inequality under the RIP condition. Cai and Zhang \cite
{CaiZhang,CaiZhang2,CaiZhang3} established the sharp bounds for the
RIP conditions that guarantee accurate recovery of low-rank matrices.

However, the RIP framework is not suitable for the ROP model considered
in the present paper. The following lemma is proved in the supplementary material~\cite{Supplement}.

\begin{Lemma}\label{lmRIPfail}
Suppose $\mathcal{X}\dvtx  \mathbb{R}^{p_1\times p_2} \to\mathbb{R}^{n}$
is ROP from the standard normal distribution. Let
\[
C_1 = \min_{A\dvtx \operatorname{rank}(A) = 1} \frac{\llVert  \mathcal{X}(A)\rrVert
_2/\sqrt{n}}{\llVert
A\rrVert  _F} \quad
\mbox{and}\quad C_2 = \max_{A\dvtx \operatorname{rank}(A) =
1}\frac{\llVert
\mathcal{X}(A)\rrVert  _2/\sqrt{n}}{\llVert  A \rrVert  _F}.
\]
Then\vspace*{1pt} for all $t>1$, $C_2/C_1 \geq\sqrt{p_1p_2/(4tn)}$ with
probability at least $1 - e^{-p_1/4} - e^{-p_2/4} - \frac{8}{n(t-1)^2}$.
\end{Lemma}

Lemma \ref{lmRIPfail} implies that at least $O(p_1p_2)$ number of
measurements are needed in order to ensure that $\mathcal{X}$
satisfies the RIP condition that guarantees the recovery of only
rank-one matrices. Since $O(p_1p_2)$ is the degree of freedom for all
matrices $A\in\mathbb{R}^{p_1\times p_2}$ and it is the number of
measurements needed to recover all $p_1\times p_2$ matrices (not just
the low-rank matrices), Lemma \ref{lmRIPfail} shows that the RIP framework is not
suitable for the ROP model.
In comparison, Theorem \ref{thgaussianl1bound} shows that if $n \geq
O(r(p_1+p_2))$, then with high probability $\mathcal{X}$ satisfies the
RUB condition of order $r$ with bounded $C_2/C_1$, which ensures the
exact recovery of all rank-$r$ matrices.

The main technical reason for the failure of RIP under the ROP model is
that RIP requires an upper bound for
%
\begin{equation}
\label{eqRIPfail} \max_{A \in\mathcal{C}}\bigl\llVert \mathcal{X}(A)\bigr
\rrVert _2^2/n = \max_{A\in
\mathcal{C}} \Biggl(\sum
_{j=1}^n \bigl(\bigl(\beta^{(j)}
\bigr)^\intercal A\gamma^{(j)} \bigr)^2 \Biggr)\Big/n,
\end{equation}
where $\mathcal{C}$ is a set containing low-rank matrices.
The right-hand side of (\ref{eqRIPfail}) involves the 4th power of
the Gaussian (or sub-Gaussian) variables $\beta^{(j)}$ and $\gamma
^{(j)}$. A much larger $n$ than the bound given in (\ref{nlowerbd})
is needed in order for the linear map $\mathcal{X}$ to satisfy the
required RIP condition, which would lead to suboptimal result.

Koltchinskii et~al. \cite{KLT} uses RIP in expectation, which is a
weaker condition than RIP. A random linear map $\mathcal{X}\dvtx  \mathbb
{R}^{p_1\times p_2} \to\mathbb{R}^{n}$ is said to satisfy RIP in
expectation of order $r$ with parameters $0<\mu<\infty$ and $0\le
\delta_r< 1$ if
\[
(1-\delta_r) \llVert A \rrVert _F^2 \leq\mu
\frac{1}{n}E\bigl\llVert \mathcal{X}(A)\bigr\rrVert _2^2
\leq(1+\delta_r)\llVert A \rrVert _F^2
\]
for all rank-$r$ matrices $A\in\mathbb{R}^{p_1\times p_2}$. This
condition was originally introduced by Koltchinskii et~al. \cite{KLT}
to prove an oracle inequality for the estimator they proposed and a
minimax lower bound. The condition is not sufficiently strong to
guarantee the exact recovery of rank-$r$ matrices in the noiseless
case. To be more specific, the bounds in Theorems 1~and~2 in \cite
{KLT} depend\vspace*{1pt} on ${\mathbf M} = \llVert  \frac{1}{n}\sum_{i=1}^n
(y_iX_i - E(y_iX_i) )\rrVert  $, which might be nonzero even in
the noiseless case. In fact, in the ROP model considered in the present
paper, we have
\begin{eqnarray*}
\frac{1}{n}E\llVert \mathcal{X}\rrVert _2^2 &= &
\frac{1}{n}\sum_{i=1}^n E \bigl(
\beta^{(i)T}A\gamma^{(i)} \bigr)^2 = E\bigl(
\beta^{\T} A\gamma\gamma^{\T} A^{\T} \beta\bigr)
\\
&= & E\operatorname{tr}\bigl(A\gamma\gamma^{\T} A^{\T}\beta
\beta^{\T}\bigr) = \operatorname{tr}\bigl(AA^{\T}\bigr)=\llVert A\rrVert
_F^2
\end{eqnarray*}
which means RIP in expectation is met for $\mu=1$ and $\delta_r=0$
for any number of measurements $n$. However, as we discussed earlier in
this section that at least $O(r(p_1+p_2))$ measurements are needed to
guarantee the model identifiability for recovery of all rank-$r$
matrices, we can see that RIP in expectation cannot ensure recovery.

Dvijotham and Fazel \cite{Dvijotham} and Oymak et al. \cite{Oymak11}
used a condition called the Spherical Section Property (SSP) which
focuses on the null space of $\mathcal{X}$. $\operatorname
{Null}(\mathcal{X})$ is
said to satisfy $\Delta$-SSP if for all $Z\in\operatorname
{Null}(\mathcal
{X})\setminus\{0\}$, $\llVert  Z\rrVert  _\ast/ \llVert  Z\rrVert  _F \geq\sqrt{\Delta}$.
Dvijotham and Fazel \cite{Dvijotham} showed that if $\mathcal{X}$
satisfies $\Delta$-SSP, $p_1\leq p_2$ and $\operatorname{rank}(A) <
\min
(3p_1/4 - \sqrt{9p_1^2/16 - p_1\Delta/4}, p_1/2 )$, the nuclear
norm minimization (\ref{eqnuclearnormminimizationnoiseless})
recovers $A$ exactly in the noiseless case. However, the SSP condition
is difficult to utilize in the ROP framework since it is hard to
characterize the matrices $Z\in\operatorname{Null}(\mathcal{X})$
when $\mathcal
{X}$ is rank-one projections.

\subsection{Gaussian noise case}\label{GaussianROPsec}


We now turn to the Gaussian noise case where $z_i\stackrel{\mathrm{i.i.d.}}{\sim}
N(0, \sigma^2)$ in (\ref{eqROP-model}). We begin by introducing a
constrained nuclear norm minimization estimator. Define two sets
%
\begin{equation}
\label{eqBl1} \mathcal{Z}_1 = \bigl\{z\dvtx  \llVert z\rrVert
_1/n\leq\sigma\bigr\} \quad\mbox{and}\quad \mathcal{Z}_2
= \bigl\{z\dvtx \bigl\llVert \mathcal{X}^\ast(z)\bigr\rrVert \leq\eta\bigr\},
\end{equation}
where $\eta= \sigma (12\sqrt{\log n}(p_1+p_2) + 6\sqrt
{2n(p_1+p_2)} )$, and let
%
\begin{equation}
\label{eqBG} \mathcal{Z}_G= \mathcal{Z}_1 \cap
\mathcal{Z}_2.
\end{equation}
Note that both $ \mathcal{Z}_1$ and $\mathcal{Z}_2$ are convex sets
and so is $\mathcal{Z}_G$. Our estimator of $A$ is given~by
%
\begin{equation}
\label{eqnuclearnormminimization2} \hat A = \mathop{\arg\min}_{M} \bigl\{\llVert M
\rrVert _\ast\dvtx  y-\mathcal{X}(M) \in \mathcal {Z}_G \bigr\}.
\end{equation}


The following theorem gives the rate of convergence for the estimator
$\hat A $ under the squared Frobenius norm loss.

\begin{Theorem}[(Upper bound)]\label{thmixedestimator}
Let $\mathcal{X}$ be ROP from the standard normal distribution and let
$z_1, \ldots, z_n \stackrel{\mathit{i.i.d.}}{\sim}N(0, \sigma^2)$. Then there
exist uniform constants $C$, $W$ and $\delta$ such that, whenever $n
\geq Cr(p_1+p_2)$, the estimator $\hat A$ given in (\ref
{eqnuclearnormminimization2}) satisfies
%
\begin{equation}
\llVert \hat A - A\rrVert _F^2 \leq W
\sigma^2 \min \biggl(\frac{r\log n
(p_1+p_2)^2}{n^2} + \frac{r(p_1+p_2)}{n}, 1 \biggr)
\end{equation}
for all rank-$r$ matrices $A$, with probability at least $1- 11/n -
3\exp(-\delta(p_1+p_2))$.
\end{Theorem}

Moreover, we have the following lower bound result for ROP.

%
\begin{Theorem}[(Lower bound)]\label{thlowerbound}
Assume that $\mathcal{X}$ is ROP from the standard normal distribution
and that $z_1, \ldots, z_n \stackrel{\mathit{i.i.d.}}{\sim}N(0, \sigma^2)$.
There exists a uniform constant $C$ such that, when $n > Cr\max(p_1,
p_2)$, with probability at least $1- 26n^{-1}$,
%
\begin{eqnarray}
&& \inf_{\hat A} \sup_{A\in\mathbb{R}^{p_1 \times p_2}\dvtx  \operatorname
{rank}(A)=r} P_z \biggl(\llVert \hat A - A\rrVert _F^2\geq\frac{\sigma^2 r(p_1 + p_2)}{32n} \biggr)
\nonumber\\[-8pt]\label{eqlowerboundP}  \\[-8pt]\nonumber
&&\qquad  \geq1 - e^{-(p_1+p_2)r/64},
\\
\label{eqlowerboundE} && \inf_{\hat A}\sup_{A\in\mathbb{R}^{p_1 \times p_2}\dvtx \operatorname
{rank}(A) = r}E_z
\llVert \hat A- A\rrVert _F^2 \geq\frac{\sigma^2r(p_1 + p_2)}{4n},
\end{eqnarray}
where $E_z$, and $P_z$ are the expectation and probability with respect
to the distribution of $z$.

When $n < r\max(p_1, p_2)$, then
%
\begin{equation}
\label{eqlowerboundinfty} \inf_{\hat A}\sup_{A\in\mathbb{R}^{p_1 \times p_2}\dvtx \operatorname
{rank}(A) =
r}E_z
\llVert \hat A - A\rrVert _F^2 =\infty.
\end{equation}
\end{Theorem}

Comparing Theorems \ref{thmixedestimator}~and~\ref{thlowerbound}, our proposed estimator is rate optimal in the
Gaussian noise case when $n \gtrsim\log n(p_1 + p_2)$ [which is
equivalent to $n\gtrsim(p_1+p_2)\log(p_1+p_2)$] or $n \sim
r(p_1+p_2)$. Since $n \gtrsim r(p_1 + p_2)$, this condition is also
implied by $r \gtrsim\log(p_1 + p_2)$. Theorem \ref{thlowerbound}
also shows that no method can recover matrices of rank $r$ consistently
if the number of measurements $n$ is smaller than $r\max(p_1, p_2)$.

The result in Theorem \ref{thmixedestimator} can also be extended to
the more general case where the matrix of interest $A$ is only
approximately low-rank. Let $A = A_{\max(r)} + A_{-\max(r)}$.

\begin{Proposition}\label{approx-low-rankprop}
Under the assumptions of Theorem \ref{thmixedestimator}, there exist
uniform constants $C$, $W_1$, $W_2$ and $\delta$ such that, whenever
$n \geq Cr(p_1+p_2)$, the estimator $\hat A$ given in (\ref
{eqnuclearnormminimization2}) satisfies
%
\begin{eqnarray}\label{eqapproxlowrankinequality}
\llVert \hat A - A\rrVert _F^2 &\leq&
W_1\sigma^2 \min \biggl(\frac{r\log n(p_1+p_2)^2}{n^2} +
\frac{r(p_1+p_2)}{n}, 1 \biggr)
\nonumber\\[-8pt]\\[-8pt]\nonumber
&&{} + W_2\frac{\llVert A_{-\max(r)}\rrVert  _\ast^2}{r}
\end{eqnarray}
for all matrices $A\in\mathbb{R}^{p_1 \times p_2}$, with probability
at least $1- 11/n - 3\exp(-\delta(p_1+p_2))$.
\end{Proposition}

If the matrix $A$ is approximately of rank $r$, then $\llVert  A_{-\max(r)}\rrVert
_\ast$ is small, and the estimator $\hat A$ continues to perform well.
This result shows that the constrained nuclear norm minimization
estimator is adaptive to the rank $r$ and robust against perturbations
of small amplitude.

%
\begin{Remark}\label{rmRademacher} All the results remain true
if the Gaussian design is replaced by the Rademacher design where
entries of $\beta^{(i)}$ and $\gamma^{(i)}$ are i.i.d. $\pm1$ with
probability ${1\over2}$. More general sub-Gaussian design case will be
discussed in Section~\ref{Extensionssec}.
\end{Remark}
%

%
\begin{Remark}\label{rmcomparsion} The estimator $\hat A$ we
propose here is the minimizer of the nuclear norm under the constraint
of the intersection of two convex sets $\mathcal{Z}_1$~and~$\mathcal
{Z}_2$. Nuclear norm minimization under either one of the two
constraints, called ``$\ell_1$~constraint nuclear norm minimization''
($\mathcal{Z} = \mathcal{Z}_1$) and ``matrix Dantzig Selector''
($\mathcal{Z}=\mathcal{Z}_2$), has been studied before in various
settings \cite
{CandesOracle,RechtMatrix,CaiZhang,CaiZhang2,CaiZhang3,Chen}. Our
analysis indicates the following:
\begin{longlist}[3.]
\item[1.] The $\ell_1$ constraint minimization performs better than the
matrix Dantzig Selector for small $n$ ($n\sim r(p_1+p_2)$) when $r \ll
\log n$.
\item[2.] The matrix Dantzig Selector outperforms the $\ell_1$ constraint
minimization for large $n$ as the loss of the matrix Dantzig Selector
decays at the rate $O(n^{-1})$.
\item[3.] The proposed estimator $\hat A$ combines the advantages of the
two estimators.
\end{longlist}
See Section~\ref{Simulationssec} for a comparison of numerical
performances of the three methods.
\end{Remark}
%

\subsection{Recovery of symmetric matrices}\label{GaussianSROPsec}

For applications such as low-dimen\-sional Euclidean embedding \cite
{Trosset,RechtMatrix}, phase retrieval \cite
{CandesRIPfail,Candes-Li} and covariance matrix estimation \cite
{Chen,CaiMaWu1,CaiMaWu2}, the low-rank matrix $A$ of interest is known
to be symmetric. Examples of such matrices include distance matrices,
Gram matrices, and covariance matrices. When the matrix $A$ is known to
be symmetric, the ROP design can be further simplified by taking $\beta
^{(i)}=\gamma^{(i)}$.

Denote\vspace*{1pt} by $ \mathbb{S}^p$ the set of all $p\times p$ symmetric
matrices in $\mathbb{R}^{p\times p}$. Let $\beta^{(1)}$, $\beta^{(2)},
\ldots, \beta^{(n)}$ be independent $p$-dimensional random vectors with
i.i.d. entries generated from some distribution $\mathcal{P}$. Define
a linear map $\mathcal{X}\dvtx  \mathbb{S}^p\to\mathbb{R}^n$ by
\[
\bigl[\mathcal{X}(A)\bigr]_i = \bigl(\beta^{(i)}
\bigr)^\intercal A\beta^{(i)} ,\qquad i=1, \ldots, n.
\]
We call such a linear map $\mathcal{X}$ ``Symmetric Rank-One
Projections'' (SROP) from the distribution $\mathcal{P}$.

Suppose we observe
%
\begin{equation}
\label{eqSROP-model} y_i =\bigl(\beta^{(i)}\bigr)^\intercal
A\beta^{(i)} +z_i,\qquad i =1, \ldots, n
\end{equation}
and wish to recover the symmetric matrix $A$. As for the ROP model, in
the noiseless case we estimate $A$ under the SROP model by
%
\begin{equation}
\label{eqnuclearminisymnoiseless} A_*= \mathop{\arg\min}_{M\in\mathbb{S}^p} \bigl\{\llVert M \rrVert
_\ast\dvtx  y = \mathcal {X}(M)\bigr\}.
\end{equation}

%
\begin{Proposition}\label{prsymmetricnoiseless}
Let $\mathcal{X}$ be SROP from the standard normal distribution.
Similar to Corollary \ref{approx-low-rankprop}, there exist uniform
constants $C$ and $\delta$ such that, whenever $n\geq Crp$, the
nuclear norm minimization estimator $A_\ast$ given by (\ref
{eqnuclearminisymnoiseless}) recovers exactly all rank-$r$
symmetric matrices $A\in\mathbb{S}^{p}$ with probability at least
$1-e^{-n\delta}$.
\end{Proposition}

For the noisy case, we propose a constraint nuclear norm minimization
estimator similar to (\ref{eqnuclearnormminimization2}). Define the
linear map $\tilde{\mathcal{X}}\dvtx  \mathbb{R}^{p_1\times p_2} \to
\mathbb{R}^{\lfloor{n}/{2} \rfloor}$ by
%
\begin{equation}
\label{eqtildeX} \bigl[\tilde{\mathcal{X}}(A)\bigr]_i = \bigl[
\mathcal{X}(A)\bigr]_{2i-1} - \bigl[\mathcal {X}(A)\bigr]_{2i},
\qquad i=1,\ldots, \biggl\lfloor\frac{n}{2}\biggr\rfloor
\end{equation}
and define $\tilde y\in\mathbb{R}^{\lfloor n/2\rfloor}$ by
%
\begin{equation}
\label{eqtildey} \tilde y_i = y_{2i-1} - y_{2i},
\qquad i=1,\ldots, \biggl\lfloor\frac
{n}{2}\biggr\rfloor.
\end{equation}
Based on the definition of $\tilde{\mathcal{X}}$, the dual map
$\tilde{\mathcal{X}}^\ast\dvtx \mathbb{R}^{\lfloor{n}/{2}\rfloor
}\to\mathbb{S}^p$ is
%
\begin{equation}
\tilde{\mathcal{X}}^\ast(z) = \sum_{i=1}^{\lfloor{n}/{2}\rfloor}
z_i \bigl(\beta^{(2i-1)}\beta^{(2i-1)\intercal} -
\beta^{(2i)}\beta^{(2i)\intercal} \bigr).
\end{equation}
Let $\eta=24\sigma (\sqrt{pn} + 2p\sqrt{2 \log n} )$.
The estimator $\hat A$ of the matrix $A$ is given by
%
\begin{equation}
\label{eqhatAsymmetric} \hat A = \mathop{\arg\min}_{M \in\mathbb{S}^p} \bigl\{\llVert M
\rrVert _\ast \dvtx  \bigl\llVert y-\mathcal{X}(M) \bigr\rrVert
_1/n \leq\sigma, \bigl\llVert \tilde{\mathcal{X}}^\ast
\bigl(\tilde y - \tilde{\mathcal{X}}(M) \bigr)\bigr\rrVert \leq\eta \bigr\}.
\end{equation}

\begin{Remark}
An important property in the ROP model considered
in Section~\ref{GaussianROPsec} is that $E\mathcal{X}=0$, that is,
$E X_i=0$ for all the measurement matrices $X_i$. However, under the
SROP model $X_i=\beta^{(i)}(\beta^{(i)})^\intercal$ and so
$E\mathcal{X}\neq0$. The step of taking the pairwise differences in
(\ref{eqtildeX}) and (\ref{eqtildey}) is to ensure that $E\tilde
{\mathcal{X}} = 0$.
\end{Remark}

The following result is similar to the upper bound given in Proposition
\ref{approx-low-rankprop} for~ROP.

\begin{Proposition}\label{prsymmetricmixedestimator}
Let $\mathcal{X}$ be SROP from the standard normal distribution and
let $z_1, \ldots, z_n\stackrel{\mathit{i.i.d.}}{\sim}N(0, \sigma^2)$. There
exist constants $C, W_1, W_2$ and $\delta$ such that, whenever $n\geq
Crp$, the estimator $\hat A$ given in (\ref{eqhatAsymmetric}) satisfies
%
\begin{equation}
\label{eqthinequalitysymmetric} \llVert \hat A - A\rrVert _F^2 \leq
W_1\sigma^2 \min \biggl(\frac{rp^2\log
n}{n^2} +
\frac{rp}{n}, 1 \biggr) + W_2\frac{\llVert  A_{-\max(r)}\rrVert  _\ast^2}{r}
\end{equation}
for all matrices $A\in\mathbb{S}^{p}$, with probability at least $1 -
15/n - 5\exp(-p\delta)$.
\end{Proposition}

In addition, we also have lower bounds for SROP, which show that the
proposed estimator is rate-optimal when
$n \gtrsim p\log n$ or $n \sim rp$, and no estimator can recover a
rank-$r$ matrix consistently if the number of measurements $n<\lfloor
\frac{r}{2}\rfloor\cdot\lfloor\frac{p}{2}\rfloor$.

\begin{Proposition}[(Lower bound)]\label{prsymmetriclowbound}
Assume that $\mathcal{X}$ is SROP from the standard normal
distribution and that $z_1, \ldots, z_n\stackrel{\mathit{i.i.d.}}{\sim
}N(0,\sigma^2)$. Then there exists a uniform constant $C$ such that,
when $n > Crp$ and $p, r\geq2$, with probability at least \mbox{$1- 26n^{-1}$},
\begin{eqnarray*}
\inf_{\hat A} \sup_{A\in\mathbb{S}^{p}\dvtx  \operatorname{rank}(A)=r} P_z
\biggl(\llVert \hat A- A\rrVert _F^2 \geq
\frac{\sigma^2 rp}{192n} \biggr) &\geq&1 - e^{-pr/192},
\\
\inf_{\hat A}\sup_{A\in\mathbb{S}^{p}\dvtx \operatorname{rank}(A) =
r}E_z
\llVert \hat A - A\rrVert _F^2 &\geq&\frac{\sigma^2rp}{24n},
\end{eqnarray*}
where $\hat A$ is any estimator of $A$, $E_z, P_z$ are the expectation
and probability with respect to $z$.

When $n < \lfloor\frac{r}{2}\rfloor\cdot\lfloor\frac{p}{2}\rfloor
$ and $p, r\geq2$, then
\[
\inf_{\hat A}\sup_{A\in\mathbb{S}^{p}\dvtx \operatorname{rank}(A) =
r}E_z
\llVert \hat A - A\rrVert _F^2 =\infty.
\]
\end{Proposition}


\section{Sub-Gaussian design and sub-Gaussian noise}
\label{Extensionssec}

We have focused on the Gaussian design and Gaussian noise distribution
in Section~\ref{Gaussiansec}. These results can be further extended
to more general distributions.
In this section, we consider the case where the ROP design is from a
symmetric sub-Gaussian distribution $\mathcal{P}$ and the errors $z_i$
are also from a sub-Gaussian distribution. We say the distribution of a
random variable $Z$ is sub-Gaussian with parameter $\tau$ if
%
\begin{equation}
\label{eqsub-Guassiannoise} P\bigl( \llvert Z\rrvert \geq t\bigr) \leq2\exp
\bigl(-t^2/\bigl(2\tau^2\bigr)\bigr)\qquad\mbox{for all
}t>0.
\end{equation}
The following lemma provides a necessary and sufficient condition for
symmetric sub-Gaussian distributions.

\begin{Lemma}\label{prsub-Gaussian}
Let $\mathcal{P}$ be a symmetric distribution and let the random
variable $X\sim\mathcal{P}$. Define
%
\begin{equation}
\label{eqalphaP} \alpha_{\mathcal{P}} = \sup_{k\geq1} \biggl(
\frac
{EX^{2k}}{(2k-1)!!} \biggr)^{{1}/{2k}}.
\end{equation}
Then the distribution $\mathcal{P}$ is sub-Gaussian if and only if
$\alpha_{\mathcal{P}}$ is finite.
\end{Lemma}




For the sub-Gaussian ROP design and sub-Gaussian noise, we estimate the
low-rank matrix $A$ by the estimator $\hat A$ given in (\ref
{eqnuclearnormminimization}) with
%
\begin{eqnarray}\label{eqBGsub-Gaussian} \qquad\mathcal{Z}_G &= & \bigl\{z\dvtx \llVert z\rrVert
_1/n \leq6\tau\bigr\}
\nonumber\\[-8pt]\\[-8pt]\nonumber
&&{} \cap \bigl\{z\dvtx \bigl\llVert \mathcal{X}^\ast(z)\bigr\rrVert \leq6
\alpha_{\mathcal
{P}}^2\tau \bigl(\sqrt{6n(p_1+p_2)}
+ 2\sqrt{\log n}(p_1+p_2) \bigr) \bigr\},
\end{eqnarray}
where $\alpha_{\mathcal{P}}$ is given in (\ref{eqalphaP}).

%
\begin{Theorem}\label{thsubgaussiannoise}
Suppose $\mathcal{X}\dvtx \mathbb{R}^{p_1\times p_2}\to\mathbb{R}^n$ is
ROP from a symmetric and variance 1 sub-Gaussian distribution $\mathcal
{P}$. Assume that $z_i$ are i.i.d. sub-\break Gaussian with parameter $\tau$
and $\hat A$ is given by (\ref{eqnuclearnormminimization}) with
$\mathcal{Z} = \mathcal{Z}_G$ defined in (\ref
{eqBGsub-Gaussian}). Then there exist constants $C, W_1, W_2, \delta
$ which only depend on $\mathcal{P}$, such that if $n\geq Cr(p_1 +
p_2)$, we have
%
\begin{eqnarray}\label{eqsubgaussiannoiseinequality}
\llVert \hat A - A\rrVert _F^2 &\leq& W_1\tau^2\min \biggl(\frac{r \log
n(p_1+p_2)^2}{n^2} +
\frac{r(p_1+p_2)}{n}, 1 \biggr)
\nonumber\\[-8pt]\\[-8pt]\nonumber
&&{} + W_2\frac{\llVert
A_{-\max(r)}\rrVert  _\ast^2}{r}
\end{eqnarray}
with probability at least $1 - 2/n - 5e^{-\delta(p_1+p_2)}$.
\end{Theorem}

%

An exact recovery result in the noiseless case for the sub-Gaussian
design follows directly from Theorem \ref{thsubgaussiannoise}. If
$z=0$, then, with high probability, all rank-$r$ matrices $A$ can be
recovered exactly via the constrained nuclear minimization (\ref
{eqnuclearnormminimizationnoiseless}) whenever $n\geq C_\mathcal
{P}r(p_1+p_2)$ for some constant $C_\mathcal{P} >0$.

%
\begin{Remark}For the SROP model considered in Section~\ref
{GaussianSROPsec}, we can similarly extend the results to the case of
sub-Gaussian design and sub-Gaussian noise. Suppose $\mathcal{X}$ is
SROP from a symmetric variance 1 sub-Gaussian distribution $\mathcal
{P}$ (other than the Rademacher $\pm$1 distribution) and $z$ satisfies
(\ref{eqsub-Guassiannoise}). Define the estimator of the low-rank
matrix $A$ by
%
\begin{equation}
\label{eqhatASROPsubGaussian} \quad\hat A = \mathop {\arg\min} _{M\in\mathbb{S}^p} \bigl\{\llVert M
\rrVert _\ast\dvtx  \bigl\llVert y-\mathcal{X}(M)\bigr\rrVert
_1/n\leq 6\tau, \bigl\llVert \tilde{\mathcal{X}}^\ast
\bigl(\tilde{y} - \tilde{\mathcal {X}}(M)\bigr)\bigr\rrVert \leq\eta \bigr\},
\end{equation}
where $\eta= C_{\mathcal{P}} (\sqrt{np} + \sqrt{\log n}p
)$ with $C_{\mathcal{P}}$ some constant depending on $\mathcal{P}$.

\begin{Proposition}\label{prexceptRademacher} Suppose $\mathcal
{X}\dvtx \mathbb{R}^{p\times p}\to\mathbb{R}^n$ is SROP from a symmetric
sub-Gaussian distribution $\mathcal{P}$ with variance 1. Also, assume
that $\operatorname{Var}(\mathcal{P}^2) > 0$ [i.e.,~$\operatorname
{Var}(w^2) > 0$ where $w \sim
\mathcal{P}$]. Let $\hat A$ be given by (\ref
{eqhatASROPsubGaussian}). Then there exist constants $ C,
C_{\mathcal{P}}, W_1, W_2$ and $\delta$ which only depend on
$\mathcal{P}$, such that for $n \geq Crp$,
%
\begin{equation}
\label{eqapproxsyslowrankinequality} \llVert \hat A - A\rrVert _F^2 \leq
W_1 \tau^2 \min \biggl(\frac{r p^2\log n}{n^2} +
\frac
{rp}{n}, 1 \biggr) + W_2 \frac{\llVert  A_{-\max(r)}\rrVert  _\ast^2}{r}
\end{equation}
with probability at least $1 - 2/n - 5e^{-\delta p}$.
\end{Proposition}

By restricting $\operatorname{Var}(\mathcal{P}^2) >0$, Rademacher
$\pm1$ is the
only symmetric and variance 1 distribution that has been excluded. The
reason why the\break Rademacher $\pm1$ distribution is an exception for the
SROP design is as follows. If $\beta^{(i)}$ are i.i.d. Rademacher $\pm
1$ distributed, then
\[
\bigl[\mathcal{X}(A)\bigr]_i = \bigl(\beta^{(i)}
\bigr)^\intercal A\beta^{(i)} = \sum_{j=1}^p
a_{jj} + \sum_{j\neq k} \beta_j^{(i)}
\beta_k^{(i)} a_{jk},\qquad i = 1,\ldots, n.
\]
So the only information contained in $\mathcal{X}(A)$ about
$\operatorname{diag}(A)$ is $\operatorname{trace}(A)$, which makes it impossible to
recover the whole matrix $A$.
\end{Remark}
\section{Application to estimation of spiked covariance matrix}\label{Applicationssec}

In this section, we consider an interesting application of the methods
and results developed in the previous sections to estimation of a
spiked covariance matrix based on one-dimensional projections. As
mentioned in the \hyperref[Introductionsec]{Introduction}, spiked covariance matrix model has been
used in a wide range of applications and it has been well studied in
the context of PCA based on i.i.d. data where one observes i.i.d.
\mbox{$p$-}dimensional random vectors $X^{(1)}, \ldots, X^{(n)}$ with mean 0
and covariance matrix $\Sigma$, where $\Sigma= I_p + \Sigma_0$ and
$\Sigma_0$ being low-rank. See, for example, \cite
{Johnstone,Birnbaum,CaiMaWu1,CaiMaWu2}. Here, we consider estimation of
$\Sigma_0$ (or equivalently $\Sigma$) based only on one-dimensional
random projections of $X^{(i)}$. More specifically, suppose that the
random vectors $X^{(1)}, \ldots, X^{(n)}$ are not directly observable
and instead we observe
%
\begin{equation}
\label{eqcovarianceobs} \xi_{i} = \bigl\langle\beta^{(i)},
X^{(i)} \bigr\rangle= \sum_{j=1}^p
\beta ^{(i)}_j X^{(i)}_j,\qquad i=1,
\ldots, n,
\end{equation}
where $\beta^{(i)}\stackrel{\mathrm{i.i.d.}}{\sim} N(0, I_p)$. The goal is to
recover $\Sigma_0$ from the projections $\{\xi_{i}, i= 1,\ldots,
n\}$.

Let $y=(y_1, \ldots, y_n)^\intercal$ with $y_i = \xi_{i}^2 - \beta
^{(i)\intercal}\beta^{(i)}$.
Note that
\[
E\bigl(\xi^2\mid \beta\bigr) =E \biggl(\sum
_{i,j} \beta_i\beta_j
X_iX_j\Big|\beta \biggr) = \sum
_{i,j} \beta_i\beta_j
\sigma_{i,j} = \beta ^\intercal\Sigma\beta
\]
and so $E(\xi^2 - \beta^\intercal\beta\mid \beta) =\beta^\intercal
\Sigma_0\beta$.
Define a linear map $\mathcal{X}\dvtx  \mathbb{S}^p \to\mathbb{R}^{n}$ by
%
\begin{equation}
\label{eqcovariancedefineMA} \bigl[\mathcal{X}(A)\bigr]_i = \beta^{(i)\intercal}A
\beta^{(i)}.
\end{equation}
Then $y$ can be formally written as
%
\begin{equation}
\label{eqyandz} y= \mathcal{X}(\Sigma_0) + z,
\end{equation}
where $z=y-\mathcal{X}(\Sigma_0)$. We define the corresponding
$\tilde{\mathcal{X}}$ and~$\tilde y$ as in (\ref{eqtildeX}) and~(\ref{eqtildey}), respectively, and apply the constraint nuclear norm
minimization to recover the low-rank matrix $\Sigma_0$ by
%
\begin{equation}
\label{eqminimizationcovariance} \hat\Sigma_0 = \mathop{\arg\min}_{M}
\bigl\{\llVert M \rrVert _\ast\dvtx  \bigl\llVert y - \mathcal {X}(M)\bigr
\rrVert \leq\eta_1, \bigl\llVert \tilde{\mathcal{X}}^*\bigl(\tilde y
- \tilde {\mathcal{X}}(M)\bigr)\bigr\rrVert \leq\eta_2 \bigr\}.
\end{equation}
The tuning parameters $\eta_1$ and $\eta_2$ are chosen as
%
\begin{equation}
\label{eqcovariancelambda} \eta_1 = c_1\sum
_{i=1}^n\xi_{i}^2 \quad
\mbox{and}\quad \eta_2 = 24c_2\sqrt{p\sum
_{i=1}^n \xi_i^4}
+ 48c_3p\log n\max_{1\leq i\leq n} \xi_{i}^2,
\end{equation}
where $c_1>\sqrt{2}$, $c_2, c_3>1$ are constants.\vadjust{\goodbreak} 


We have the following result on the estimator (\ref
{eqminimizationcovariance}) for spiked covariance matrix estimation.

\begin{Theorem}\label{thcovariance}
Suppose $n\geq3$, we observe $\xi_{i}, i=1,\ldots, n$, as in (\ref
{eqcovarianceobs}), where $\beta^{(i)} \stackrel{\mathit{i.i.d.}}{\sim}
N(0,I_p)$ and $X^{(1)}, \ldots, X^{(n)} \stackrel{\mathit{i.i.d.}}{\sim}
N(0,\Sigma)$ with $\Sigma= I_p + \Sigma_0$ and $\Sigma_0$ positive
semidefinite and $\operatorname{rank}(\Sigma_0)\leq r$. Let $\hat\Sigma_0$ be given
by (\ref{eqminimizationcovariance}). Then there exist uniform
constants $C$, $D$, $\delta$ such that when $n\geq Drp$,
%
\begin{eqnarray}\label{eqcovarianceineq}
&& \llVert \hat\Sigma_0 - \Sigma_0\rrVert
_F^2
\nonumber\\[-8pt]\\[-8pt]\nonumber
&&\qquad \leq C\min \biggl(\frac{rp}{n}\llVert \Sigma
\rrVert _\ast^2 + \frac{r p^2\log^4 n}{n^2} \bigl(\llVert \Sigma
\rrVert _\ast^2 + \log^2 n \llVert \Sigma\rrVert
^2\bigr), \llVert \Sigma\rrVert _\ast^2 \biggr)
\end{eqnarray}
with probability at least $1-O(1/n) - 4\exp(-p\delta) - \frac
{2}{\sqrt{2\pi\log n}}$.
\end{Theorem}

%
\begin{Remark}
We have focused estimation of spiked covariance
matrices on the setting where the random vectors $X^{(i)}$ are
Gaussian. Similar to the discussion in Section~\ref{Extensionssec},
the results given here can be extended to more general distributions
under certain moment conditions.
\end{Remark}
%

%
\begin{Remark}
The problem considered in this section is related
to the so-called covariance sketching problem considered in Dasarathy
et~al. \cite{Nowak1}. In covariance sketching, the goal is to estimate
the covariance matrix of high-dimensional random vectors $X^{(1)},
\ldots, X^{(n)}$ based on the low-dimensional projections
\[
y^{(i)}=Q X^{(i)},\qquad i=1, \ldots, n,
\]
where $Q$ is a fixed $m\times p$ projection matrix with $m<p$. The main
differences between the two settings are that the projection matrix in
covariance sketch is the same for all $X^{(i)}$ and the dimension $m$
is still relatively large with $m\ge C \sqrt{p}\log^3 p$ for some
$C>0$. In our setting, $m=1$ and $Q$ is random and varies with $i$. The
techniques for solving the two problems are very different. Comparing
to \cite{Nowak1}, the results in this section indicate that there is a
significant advantage to have different random projections for
different random vectors $X^{(i)}$ as opposed to having the same
projection for all $X^{(i)}$.
\end{Remark}
%


\section{Simulation results}\label{Simulationssec}

The constrained nuclear norm minimization methods can be efficiently
implemented. The estimator $\hat A$ proposed in Section~\ref
{GaussianROPsec} can be implemented by the following convex programming:
%
\begin{eqnarray}
\mbox{minimize}&\qquad&\operatorname{Tr}(B_1) + \operatorname{Tr}(B_2)\nonumber
\\
\mbox{subject to}&\qquad&
\left[\matrix{ B_1 & A
\cr
A^T & B_2}\right]\succeq0,\qquad\bigl\llVert y-\mathcal{X}(A)
\bigr\rrVert _1 \leq\lambda _1,
\\
&\qquad& \bigl\llVert
\mathcal{X}^\ast\bigl(y-\mathcal{X}(A)\bigr)\bigr\rrVert \leq
\lambda_2,\nonumber
\end{eqnarray}
with optimization variables $B_1\in\mathbb{S}^{p_1}, B_2\in\mathbb
{S}^{p_2}$, $A\in\mathbb{R}^{p_1\times p_2}$.
We use the CVX package \cite{CVX1,CVX2} to implement the proposed
procedures. In this section, a simulation study is carried out to
investigate the numerical performance of the proposed procedures for
low-rank matrix recovery in various settings.

\begin{figure}

\includegraphics{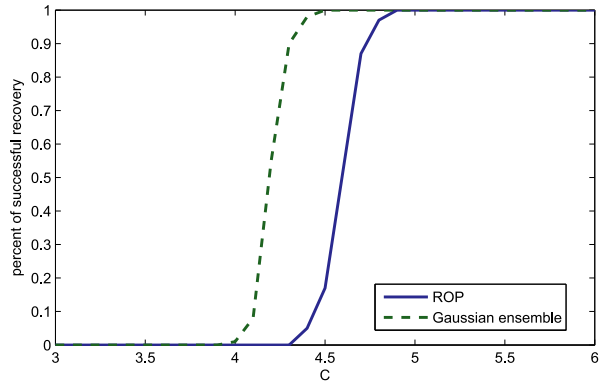}

\caption{Rates of successful recovery for the ROP and Gaussian
ensemble with $p_1 = p_2 = 100$, $r = 5$, and $n = Cr\max(p_1, p_2)$
for $C$ ranging from 3 to 6.}\label{figRORM100}
\end{figure}

We begin with the noiseless case. In this setting, Theorem \ref
{thgaussianl1bound} and Corollary~\ref{rmGaussiannoiseless} show
that the nuclear norm minimization recovers a rank $r$ matrix exactly whenever
%
\begin{equation}
\label{C} n \geq Cr\max(p_1, p_2).
\end{equation}
A similar result holds for the Gaussian ensemble \cite{CandesOracle}.
However, the minimum constant $C$ that guarantees the exact recovery
with high probability is not specified in either case. It is of
practical interest to find the minimum constant $C$. For this purpose,
we randomly generate $p_1\times p_2$ rank-$r$ matrices $A$ as $ A =
X^\intercal Y$, where $X\in\mathbb{R}^{r\times p_1}$, $Y\in\mathbb
{R}^{r\times p_2}$ are i.i.d. Gaussian matrices. We compare ROP from
the standard Gaussian distribution and the Gaussian ensemble, with the
number of measurements $n = Cr\max(p_1, p_2)$ from a range of values
of $C$ using the constrained nuclear norm minimization (\ref
{eqnuclearnormminimizationnoiseless}).
A recovery is considered successful if $\llVert  \hat A - A\rrVert  _F/\llVert  A \rrVert  _F\leq
10^{-4}$. Figure~\ref{figRORM100} shows the rate of successful
recovery when $p_1 = p_2= 100$ and $r = 5$.

%

The numerical results show that for ROP from the Gaussian distribution,
the minimum constant $C$ to ensure exact recovery with high probability
is slightly less than 5 in the small scale problems ($p_1, p_2 \leq
100$) we tested. The corresponding minimum constant $C$ for the
Gaussian ensemble is about $4.5$. Matrix completion requires much
larger number of measurements. Based on the theoretical analyses given
in \cite{CandesRecht,RechtImproved}, the required number of
measurements for matrix completion is $O(\mu r(p_1+p_2)\log
^2(p_1+p_2))$, where $\mu\geq1$ is some coherence constant describing
the ``spikedness'' of the matrix $A$. Hence, for matrix completion, the
factor $C$ in
(\ref{C}) needs to grow with the dimensions $p_1$ and $p_2$ and it\vadjust{\goodbreak}
requires $C \gtrsim\mu\log^2(p_1+p_2)$, which is much larger than
what is needed for the ROP or Gaussian ensemble. The required storage
space for the Gaussian ensemble is much greater than that for the ROP.
In order to ensure accurate recovery of $p\times p$ matrices of rank
$r$, one needs at least $4.5rp^3$ bytes of space to store the
measurement matrices, which could be prohibitively large for the
recovery of high-dimensional matrices. In contrast, the storage space
for the projection vectors in ROP is only $10rp^2$ bytes, which is far
smaller than what is required by the Gaussian ensemble in the high-dimensional case.

\begin{figure}

\includegraphics{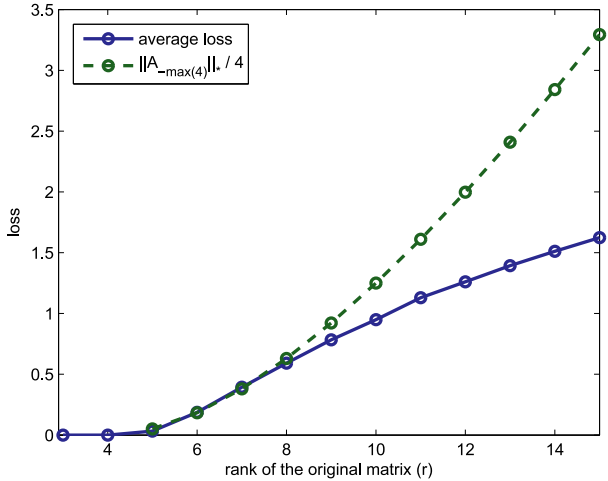}

\caption{Recovery accuracy (solid line) for approximately
low-rank matrices with different values of~$r$, where $p_1 = p_2 =
100$, $n = 2000$, $\sigma(A) = (1,1/\sqrt{2}, \ldots, 1/\sqrt{r})$.
The dashed line is the theoretical upper bound.}\label{figapproxlowrank}
\end{figure}

We then consider the recovery of approximately low-rank matrices to
investigate the robustness of the method against small perturbations.
To this end, we randomly draw $100\times100$ matrix $A$ as $A=U\cdot
\operatorname{diag}(1, 2^{-1/2}, \ldots, r^{-1/2})\cdot V^\intercal
$, where $U\in\mathbb{R}^{100\times r}$ and $V\in\mathbb
{R}^{100\times r}$ are random matrices with orthonormal columns. We
then observe $n=2000$ random rank-one projections with the measurement
vectors being i.i.d. Gaussian. Based on the observations, the nuclear
minimization procedure (\ref{eqnuclearnormminimizationnoiseless})
is applied to estimate $A$. The results for different values of $r$ are
shown in Figure~\ref{figapproxlowrank}. It can be seen from the
plot that in this setting one can exactly recover a matrix of rank at
most 4 with 2000 measurements. However, when the rank $r$ of the true
matrix $A$ exceeds 4, the estimate is still stable. The theoretical
result in Proposition \ref{approx-low-rankprop} bounds the loss (solid
line) at $O(\llVert  A_{-\max(4)}\rrVert  ^2_\ast/4)$ (shown in the dashed
line) with high probability, which corresponds to Figure~\ref
{figapproxlowrank}.
%

%
\begin{figure}[b]

\includegraphics{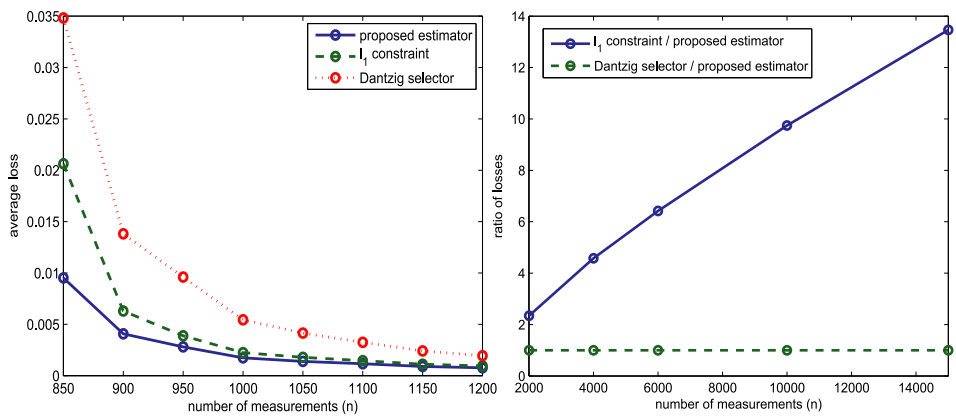}

\caption{Left panel: Comparison of the proposed estimator with $\hat A{}^{\ell_1}$ and $\hat A{}^{\mathrm{DS}}$ for $p_1 = p_2 = 50$, $r = 4$, $\sigma=
0.01$, and $n$ ranging from 850 to 1200.
Right panel: Ratio of the squared Frobenius norm loss of $\hat A{}^{\ell
_1}$ to that of the proposed estimator for $p_1 = p_2 = 50$, $r = 4$,
and $n$ varying from 2000 to~15,000.}\label{fig50comparison}
\end{figure}

We now turn to the noisy case. The low-rank matrices $A$ are generated
by $A = X^\intercal Y$, where $X\in\mathbb{R}^{r\times p_1}$ and
$Y\in\mathbb{R}^{r\times p_2}$ are i.i.d. Gaussian matrices. The ROP
$\mathcal{X}$ is from the standard Gaussian distribution and the noise
vector $z\sim N_n(0, \sigma^2)$. Based on $(\mathcal{X}, y)$ with $y
= \mathcal{X}(A)+z$, we compare our proposed estimator $\hat A$
with the $\ell_1$ constraint minimization estimator $\hat A{}^{\ell_1}$
\cite{Chen} and the matrix Dantzig Selector $\hat A{}^{\mathrm{DS}}$~\cite{CandesOracle}, where
\begin{eqnarray*}
\hat A &=& \mathop{\arg\min}_{M}\bigl\{ \llVert M \rrVert
_\ast\dvtx  y - \mathcal{X}(M) \in\mathcal {Z}_1 \cap
\mathcal{Z}_2\bigr\},
\\
 \hat A{}^{\ell_1} &=& \mathop{\arg\min}_{M}\bigl\{ \llVert M
\rrVert _\ast\dvtx  y - \mathcal{X}(M) \in\mathcal{Z}_1\bigr\},
\\
 \hat A{}^{\mathrm{DS}} &=& \mathop{\arg\min}_{M}\bigl\{ \llVert M
\rrVert _\ast\dvtx  y - \mathcal {X}(M)\in \mathcal{Z}_2\bigr\},
\end{eqnarray*}
with $\mathcal{Z}_1 = \{z\dvtx  \llVert  z\rrVert  _1/n \leq\sigma\} $ and $\mathcal
{Z}_2 = \{z\dvtx  \llVert  \mathcal{X}(z)\rrVert   \leq\sigma(\sqrt{\log n}(p_1 + p_2)
+\break  \sqrt{n(p_1 + p_2)})\}$. Note that $\hat A{}^{\ell_1}$ is similar to
the estimator proposed in Chen et~al. \cite{Chen}, except their
estimator is for symmetric matrices under the SROP but ours is for
general low-rank matrices under the ROP.
Figure~\ref{fig50comparison} compares the performance of the three
estimators. It can be seen from the left panel that for small $n$,
$\ell_1$~constrained minimization outperforms the matrix Dantzig
Selector, while our estimator outperforms both $\hat A{}^{\ell_1}$ and
$\hat A{}^{\mathrm{DS}}$. When $n$ is large, our estimator and $\hat A{}^{\mathrm{DS}}$ are
essentially the same and both outperforms $\hat A{}^{\ell_1}$. The right
panel of Figure~\ref{fig50comparison} plots the ratio of the squared
Frobenius norm loss of $\hat A{}^{\ell_1}$ to that of our estimator. The
ratio increases with $n$. These numerical results are consistent with
the observations made in Remark \ref{rmcomparsion}.

\begin{figure}

\includegraphics{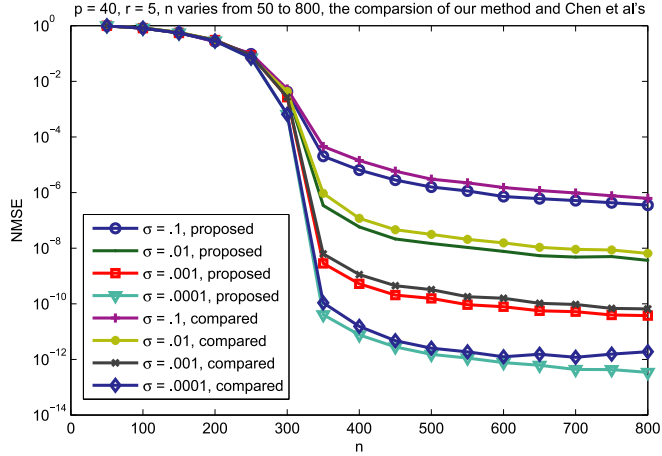}

\caption{Comparison of the proposed estimator $\hat A$ with the $\hat A{}^{\ell_1}$. Here $p = 40$, $r = 5$, $\sigma= 0.1, 0.01, 0.001,
0.0001$ and $n$ ranges from 50 to 800.}\label{fig40comparisonchen}
\end{figure}

We now turn to the recovery of symmetric low-rank matrices under the
SROP model (\ref{eqSROP-model}).
Let $ \mathcal{X}$ be SROP from the standard normal distribution. We
consider the setting where $p = 40$, $n$ varies from 50 to 600, $z_i
\sim\sigma\cdot\mathcal{U}[-1, 1]$ with $\sigma={}$0.1, 0.01, 0.001 or
0.0001, and $A$ is randomly generated as rank-5 matrix by the same
procedure discussed above. The setting is identical to the one
considered in Section~5.1 of \cite{Chen}. Although we cannot exactly
repeat the simulation study in \cite{Chen} as they did not specify the
choice of the tuning parameter, we can implement both our procedure
\begin{eqnarray*}
\hat A &=& \mathop{\arg\min}_{M} \biggl\{ \llVert M \rrVert
_\ast\dvtx   \bigl\llVert y - \mathcal{X}(M)\bigr\rrVert _1
\leq\frac{n\sigma}{2},
\\
&&\hspace*{39pt}{}\bigl\llVert \tilde{\mathcal{X}}^\ast\bigl(\tilde y - \tilde{
\mathcal{X}}(M) \bigr)\bigr\rrVert \leq\frac{\sigma(\sqrt{\log n}p + \sqrt{np})}{3} \biggr\}
\end{eqnarray*}
and the estimator $\hat A{}^{\ell_1}$ with only the $\ell_1$ constraint
which was proposed by Chen et~al. \cite{Chen}
\[
\hat A{}^{\ell_1} = \mathop{\arg\min}_{M} \biggl\{ \llVert M
\rrVert _\ast\dvtx  \bigl\llVert y - \mathcal{X}(M)\bigr\rrVert
_1 \leq\frac{n\sigma}{2} \biggr\}.
\]
The\vspace*{1pt} results are given in Figure~\ref{fig40comparisonchen}. It can
be seen that our estimator $\hat A$ outperforms the estimator $\hat A{}^{\ell_1}$.
%

\subsection{Data driven selection of tuning parameters}

We have so far considered the estimators
%
\begin{eqnarray}
\label{eqmethodROP} \hat A &=& \mathop{\arg\min}_{B}\bigl\{ \llVert B\rrVert
_\ast\dvtx  \bigl\llVert y-\mathcal{X}(B)\bigr\rrVert _1/n\leq
\lambda, \bigl\llVert \mathcal{X}^\ast\bigl(y-\mathcal{X}(B)\bigr)\bigr
\rrVert \leq\eta\bigr\},
\\
\label{eqmethodSROP} \hat A &=& \mathop{\arg\min}_{M} \bigl\{\llVert M
\rrVert _\ast\dvtx  \bigl\llVert y-\mathcal{X}(M)\bigr\rrVert
_1/n\leq \lambda, \bigl\llVert \tilde{\mathcal{X}}^\ast
\bigl(\tilde y - \tilde{\mathcal {X}}(M)\bigr)\bigr\rrVert \leq\eta\bigr\}
\end{eqnarray}
for the ROP and SROP, respectively.
The theoretical choice of the tuning parameters $\lambda$ and $\eta$
depends on the knowledge of the error distribution such as the
variance. In real applications, such information may not be available
and/or the\vadjust{\goodbreak} theoretical choice may not be the best. It is thus desirable
to have a data driven choice of the tuning parameters. We now introduce
a practical method for selecting the tuning parameters using $K$-fold
cross-validation.

Let $(\mathcal{X}, y) = \{(X_i, y_i), i=1,\ldots, n\}$ be the
observed sample and let $T$ be a grid of positive real values.
For each $t\in T$, set
%
\begin{eqnarray}
\label{eqlambdaeta}
(\lambda, \eta) &=& \bigl(\lambda(t), \eta(t)\bigr)
\nonumber\\[-4pt]\\[-12pt]\nonumber
&=& \cases{
\bigl(t, t \bigl(\sqrt{\log n}(p_1 + p_2) +
\sqrt{n(p_1 + p_2)} \bigr) \bigr), &\quad for ROP;
\vspace*{5pt}\cr
\bigl(t, t (\sqrt{\log n}p + \sqrt{np} ) \bigr), &\quad for SROP.}
\end{eqnarray}
Randomly split the $n$ samples $(X_i, y_i), i=1,\ldots, n$ into two
groups of sizes $n_1 \sim\frac{(K-1)n}{K}$ and $n_2 \sim\frac
{n}{K}$ for $I$ times. Denote by $J_1^{i}, J_2^{i} \subseteq\{1,
\ldots, n\}$ the index sets for Groups 1~and~2, respectively, for
the $i$th split. Apply our procedure [(\ref{eqmethodROP}) for ROP
and (\ref{eqmethodSROP}) for SROP, resp.] to the sub-samples
in Group~1 with the tuning parameters $(\lambda(t),\eta(t))$ and
denote the estimators by $\hat A{}^{i}(t)$, $i=1,\ldots, I$. Evaluate
the prediction error of $\hat A{}^{i}(t)$ over the subsample in Group~2
and set
\[
\hat R(t) = \sum_{i=1}^I \sum
_{j\in J_2^i} \bigl\llvert y_j - \bigl\langle
A^{i}(t), X_j\bigr\rangle\bigr\rrvert ^2,
\qquad t\in T.
\]
We select
\[
t_* = \mathop{\arg\min}_{T} \hat R(t)
\]
and choose the tuning parameters $(\lambda(t_*), \eta(t_*))$ as in
(\ref{eqlambdaeta}) with $t=t_*$ and the final estimator $\hat A$
based on (\ref{eqmethodROP}) or (\ref{eqmethodSROP}) with the
chosen tuning parameters.

%
\begin{figure}

\includegraphics{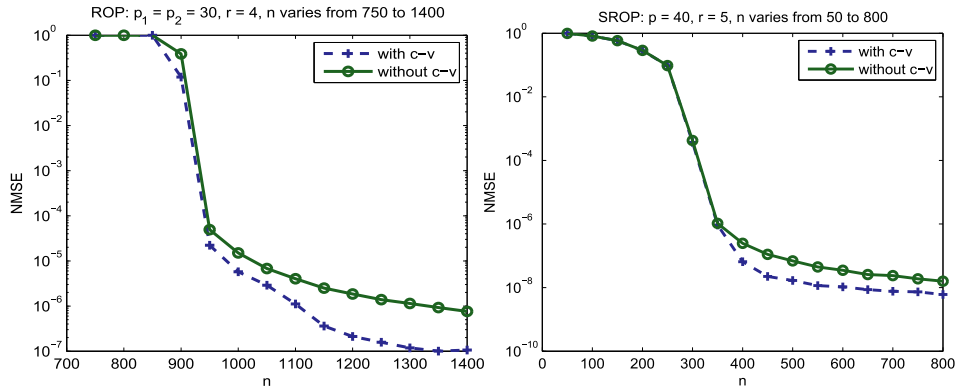}

\caption{Comparison of the performance with cross validation and
without cross-validation in both ROP and SROP. Left panel: ROP, $p_1 =
p_2 = 30$, $r=4$, $n$ varies from 750 to 1400. Right panel: SROP,
$p=40$, $r = 5$, $n$ varies from 50 to 800.}\label{figcv}
\end{figure}

We compare the numerical result by 5-fold cross-validation with the
result based on the known $\sigma$ by simulation in Figure~\ref
{figcv}. Both the ROP and SROP are\vadjust{\goodbreak} considered. It can be seen that the
estimator with the tuning parameters chosen through \mbox{5-}fold
cross-validation has the same performance as or outperforms the one
with the theoretical choice of the tuning parameters.

\subsection{Image compression}

Since a two-dimensional image can be considered as a matrix, one
approach to image compression is by using low-rank matrix approximation
via the singular value decomposition. See, for example, \cite
{Andrews,RechtMatrix,Wakinimage}. Here, we use an image recovery
example to further illustrate the nuclear norm minimization method
under the ROP model.

\begin{figure}[b]

\includegraphics{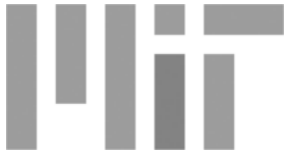}

\caption{Original grayscale MIT logo.}\label{figmit-gray}
\end{figure}

For a grayscale image, let $A=(a_{i,j}) \in\mathbb{R}^{m\times n}$ be
the intensity matrix associated with the image, where $a_{ij}$ is the
grayscale intensity of the $(i, j)$ pixel. When the matrix $A$ is
approximately low-rank, the ROP model and nuclear norm minimization
method can be used for image compression and recovery.
To illustrate this point, let us consider the following grayscale MIT
Logo image (Figure~\ref{figmit-gray}).
%

\begin{figure}

\includegraphics{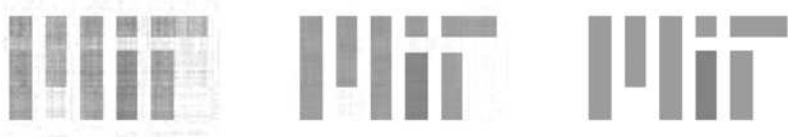}

\caption{Recovery of MIT logo based on different number of
measurements. Left: 900; Middle: 1000; Right: 1080.}\label{figmit-recovered}
\end{figure}

The matrix associated with MIT logo is of the size $50\times80$ and of
rank 6. We take rank-one random projections $\mathcal{X}(A)$ as the
observed sample, with various sample sizes. Then the constrained
nuclear norm minimization method is applied to reconstruct the original
low-rank matrix. The recovery results are shown in Figure~\ref
{figmit-recovered}. The results show that the original image can be
compressed and recovered well via the ROP model and the nuclear norm
minimization.
%

\section{Discussions}\label{Discussionssec}

This paper introduces the ROP model for the recovery of general
low-rank matrices. A constrained nuclear norm minimization method is
proposed and its theoretical and numerical properties are studied. The
proposed estimator is shown to be rate-optimal when the number of
rank-one projections $n \gtrsim\log n(p_1 + p_2)$ or $n \sim
r(p_1+p_2)$. It is also shown that the procedure is adaptive to the
rank and robust against small perturbations. The method and results are
applied to estimation of a spiked covariance matrix. It is somewhat
unexpected that it is possible to accurately recover a spiked
covariance matrix from only one-dimensional projections. An interesting
open problem is to estimate the principal components/subspace based on
the one-dimensional random projections. We leave this as future work.

In a recent paper, Chen et~al. \cite{Chen} considered quadratic
measurements for the recovery of symmetric positive definite matrices,
which is similar to the special case of SROP that we studied here. The
paper was posted on arXiv as we finish writing the present paper. They
considered the noiseless and $\ell_1$ bounded noise cases
and introduced the so-called ``RIP-$\ell_2$/$\ell_1$'' condition.
The ``RIP-$\ell_2$/$\ell_1$'' condition is similar to RUB in our
work. But these two conditions are not identical as the RIP-$\ell
_2$/$\ell_1$ condition can only be applied to symmetric low-rank
matrices as only symmetric operators are considered in the paper. In
contrast, RUB applies to all low-rank matrices.

Chen et~al. (\cite{Chen} version 4) considered $\ell_1$-bounded noise
case under the SROP model and gave an upper bound in their Theorem~3
(after a slight change of notation)
%
\begin{equation}
\label{ineqhatSigma-Sigma} \llVert \hat{\Sigma} - \Sigma\rrVert _F \leq
C_1 \frac{\llVert  \Sigma- \Sigma
_{\Omega}\rrVert  _\ast}{\sqrt{r}} + C_2\frac{\varepsilon}{n}.
\end{equation}
This result for $\ell_1$ bounded noise case is not applicable to the
i.i.d. random noise setting. When the entries of the noise term $\eta
\in\mathbb{R}^n$ are of constant order, which is the typical case for
i.i.d. noise with constant variance, one has $\llVert  \eta\rrVert  _1\sim Cn$ with
high probability. In such a case, the term $C_2\frac{\varepsilon_1}{n}$
on the right-hand side of (\ref{ineqhatSigma-Sigma}) does not even
converge to 0 as the sample size $n\to\infty$.

In comparison, the bound (\ref{eqapproxsyslowrankinequality}) in
Proposition \ref{prexceptRademacher} can be equivalently rewritten~as
%
\begin{equation}
\label{eqapproxlowrankinequalityedit} \qquad\llVert \hat A - A\rrVert _F \leq W_2
\frac{\llVert  A_{-\max(r)}\rrVert  _\ast}{\sqrt{r}} + W_1\tau\min \biggl(\frac{\sqrt{r\log n} p}{n} + \sqrt{
\frac
{rp}{n}}, 1 \biggr),
\end{equation}
where the first term $W_2\frac{\llVert  A_{-\max(r)}\rrVert  _\ast}{\sqrt{r}}$ is
of the same order as $C_1\frac{\llVert  \Sigma- \Sigma_{\Omega}\rrVert  _\ast
}{\sqrt{r}}$ in (\ref{ineqhatSigma-Sigma}) while the second term
decays to 0 as $n\to\infty$. Hence, for the recovery of rank-$r$
matrices, as the sample size $n$ increases our bound decays to 0 but
the bound (\ref{ineqhatSigma-Sigma}) given in Chen et~al. \cite
{Chen} does not.
The main reason of this phenomenon lies in the difference in the two
methods: we use nuclear norm minimization under two convex constraints
(see Remark \ref{rmcomparsion}), but Chen et~al. \cite{Chen} used only the $\ell_1$
constraint. Both theoretical results (see Remark \ref{rmcomparsion})
and numerical results (Figure~\ref{fig50comparison} in Section~\ref
{Simulationssec}) show that the additional constraint $\mathcal{Z}_2$
improves the performance of the estimator.

Moreover, the results and techniques in \cite{Chen} for symmetric
positive definite matrices are not applicable to the recovery of
general nonsymmetric matrices.
This is due to the fact that for a nonsymmetric square matrix
$A=(a_{ij})$, the quadratic measurements $(\beta^{(i)})^\intercal A
\beta^{(i)}$ satisfy
\[
\bigl(\beta^{(i)}\bigr)^\intercal A \beta^{(i)} = \bigl(
\beta^{(i)}\bigr)^\intercal A^s \beta^{(i)},
\]
where $A^s={1\over2}(A+A^\intercal)$. Hence,\vspace*{1pt} for a nonsymmetric
matrix $A$, only its symmetrized version $A^s$ can be possibly
identified and estimated based on the quadratic measurements, the
matrix $A$ itself is neither identifiable nor estimable.



\section{Proofs}\label{Proofssec}

We prove the main results in this section. We begin by collecting a few
important technical lemmas that will be used in the proofs of the main
results. The proofs of some of these technical lemmas are involved and
are postponed to the supplementary material \cite{Supplement}.

\subsection{Technical tools}

Lemmas \ref{lmmgf} and \ref{lmmean} below are used for deriving the
RUB condition (see Definition \ref{dfl1bound}) from the ROP design.

\begin{Lemma}\label{lmmgf}
Suppose $A\in\mathbb{R}^{p_1 \times p_2}$ is a fixed matrix and
$\mathcal{X}$ is ROP from a symmetric sub-Gaussian distribution
$\mathcal{P}$, that is,
\[
\bigl[\mathcal{X}(A)\bigr]_j = \beta^{(j)T} A
\gamma^{(j)},\qquad j=1,\ldots, n,
\]
where $\beta^{(j)} = (\beta_1^{(j)}, \ldots, \beta_{p_1}^{(j)})^T,
\gamma^{(j)} = (\gamma_1^{(j)},\ldots, \gamma_{p_2}^{(j)})^T$ are
random vectors with entries i.i.d. generated from $\mathcal{P}$.
Then for $\delta>0$, we have
\[
\biggl(\frac{1}{3\alpha^4_{\mathcal{P}}}-2\alpha^2_{\mathcal
{P}}\delta-
\alpha^2_{\mathcal{P}}\delta^2 \biggr)\llVert A \rrVert
_F \leq\bigl\llVert \mathcal{X}(A)\bigr\rrVert _1/n
\leq \bigl(1 + 2\alpha_{\mathcal{P}}^2\delta +\alpha_{\mathcal{P}}^2
\delta^2 \bigr)\llVert A \rrVert _F
\]
with probability at least $1- 2\exp(-\delta^2 n)$. Here, $\alpha
_\mathcal{P}$ is defined by (\ref{eqalphaP}).
\end{Lemma}

%
\begin{Lemma}\label{lmmean}
Suppose $A\in\mathbb{R}^{p_1\times p_2}$ is a fixed matrix. $\beta
=(\beta_1,\ldots,\beta_{p_1})^T, \gamma= (\gamma_1, \ldots,
\gamma_{p_2})^T$ are random vectors such that $\beta_1,\ldots, \beta
_{p_1}, \gamma_1,\ldots, \gamma_{p_2} \stackrel{\mathit{i.i.d.}}{\sim}
\mathcal{P}$, where $\mathcal{P}$ is some symmetric variance 1
sub-Gaussian distribution, then we have
\[
\frac{ \llVert  A \rrVert  _F}{3\alpha^4_{\mathcal{P}}} \leq E\bigl\llvert \beta^T A \gamma\bigr\rrvert
\leq\llVert A \rrVert _F,
\]
where $\alpha_\mathcal{P}$ is given by (\ref{eqalphaP}).\vadjust{\goodbreak}
\end{Lemma}

Let $z\in\mathbb{R}^n$ be i.i.d. sub-Gaussian distributed. By measure
concentration theory, $\llVert  z\rrVert  _p^p/n$, $1\leq p\leq\infty$, are
essentially bounded; specifically, we have the following lemma.

\begin{Lemma}\label{lml1norm}
Suppose $z\in\mathbb{R}^{n}$ and $z_i\stackrel{\mathit{i.i.d.}}{\sim} N(0,
\sigma^2)$, we have
\begin{eqnarray*}
P\bigl(\llVert z\rrVert _1 \geq\sigma n\bigr)&\leq&\frac{9}{n},
\\
P \bigl(\llVert z\rrVert _2\geq\sigma\sqrt{n+2\sqrt{n\log n}} \bigr)
&\leq&\frac{1}{n},
\\
P\bigl(\llVert z\rrVert _\infty\geq2\sigma\sqrt{\log n}\bigr) &\leq&
\frac{1}{n\sqrt{2\pi\log n}}.
\end{eqnarray*}
More general, when $z_i$ are i.i.d. sub-Gaussian distributed such that
(\ref{eqsub-Guassiannoise}) holds, then
\begin{eqnarray*}
P \bigl(\llVert z\rrVert _1 \geq Cn \bigr) &\leq&\exp \biggl(-
\frac{n(C - 2\sqrt
{2\pi}\gamma)^2}{2\gamma^2} \biggr)\qquad\forall C > 2\sqrt{2\pi }\gamma,
\\
P\bigl(\llVert z\rrVert _2 \geq\sqrt{Cn}\bigr) &\leq&\exp \biggl(-
\frac{n(C-4\gamma
^2)^2}{8\gamma^2 C} \biggr)\qquad\forall C>4\gamma^2,
\\
P\bigl(\llVert z\rrVert _\infty \geq C\gamma\sqrt{\log n}\bigr)&\leq& 2
n^{-C^2/2-1}\qquad \forall C>0.
\end{eqnarray*}
\end{Lemma}

Lemma \ref{lmDSbound} below presents an upper bound for the spectral
norm of $\mathcal{X}(z)$ for a fixed vector $z$.

\begin{Lemma}\label{lmDSbound}
Suppose $\mathcal{X}$ is ROP from some symmetric sub-Gaussian
distribution $\mathcal{P}$ and $z\in\mathbb{R}^{n}$ is some fixed
vector, then for $C>\log7$, we have
\[
\bigl\llVert \mathcal{X}^\ast(z)\bigr\rrVert \leq3
\alpha_{\mathcal
{P}}^2 \bigl(C(p_1+p_2)
\llVert z\rrVert _\infty+ \sqrt{2C(p_1+p_2)}
\llVert z\rrVert _2 \bigr) %
\]
with probability at least $1 - 2\exp (-(C - \log
7)(p_1+p_2) )$. Here, $\alpha_\mathcal{P}$ is defined by~(\ref{eqalphaP}).
\end{Lemma}

We are now ready to prove the main results of the paper.

\subsection{Proof of Theorem \texorpdfstring{\protect\ref{thl1boundrecovery}}{2.1}}
We introduce the following two technical lemmas that will be used in
the proof of theorem.

The \emph{null space property} below is a well-known result in affine
rank minimization problem (see \cite{Oymak11}). It provides a
necessary, sufficient and easier-to-check condition for exact recovery
in the noiseless setting.

\begin{Lemma}[(Null space property)]\label{lmnullspace}
Using (\ref{eqnuclearnormminimizationnoiseless}), one can recover
all matrices $A$ of rank at most $r$ if and only if for all $R\in
\mathcal{N}(\mathcal{X})\setminus\{0\}$,
\[
\llVert R_{\max(r)}\rrVert _\ast<\llVert R_{-\max(r)}
\rrVert _\ast.
\]
\end{Lemma}

The following lemma is given in \cite{CaiZhang3}, which provides a
way to decompose the general vectors to sparse ones.

\begin{Lemma}[(Sparse representation of a polytope)]\label{lmsparsemean}
Suppose $s$ is a nonnegative integer, $v\in\mathbb{R}^p$ and $\theta
\geq0$. Then $\llVert  v\rrVert  _\infty\leq\theta, \llVert  v\rrVert  _1\leq s\theta$, if and
only if $v$ can be expressed as a weighted mean,
\[
v=\sum_{i=1}^N\lambda_i
u_i,\qquad0\leq\lambda_i\leq1,\qquad\sum
_{i=1}^N \lambda_i=1,
\]
where $u_i$ satisfies%
\begin{eqnarray}\label{eqmeanu}
&& u_i \mbox{ is $s$-sparse},\qquad
\operatorname{supp}(u_i) \subseteq \operatorname{supp} (v),
\nonumber\\[-8pt]\\[-8pt]\nonumber
&& \llVert u_i\rrVert _1=\llVert v\rrVert _1, \qquad\llVert u_i\rrVert _\infty\leq\theta.
\end{eqnarray}
\end{Lemma}

For the proof of Theorem \ref{thl1boundrecovery}, by null space
property (Lemma \ref{lmnullspace}), we only need to show for all
nonzero $R$ with $\mathcal{X}(R)=0$, we must have $\llVert  R_{\max(r)}\rrVert
_\ast< \llVert  R_{-\max(r)}\rrVert  _\ast$.

If this does not hold, suppose there exists nonzero $R$ with $\mathcal
{X}(R)=0$ and $\llVert  R_{\max(r)}\rrVert  _\ast\geq\llVert  R_{-\max(r)}\rrVert  _\ast$. We
denote $p=\min(p_1, p_2)$ and assume the singular value decomposition
of $R$ is
\[
R=\sum_{i=1}^{p} \sigma_i
u_i v_i^\T = U\operatorname{diag}(\vec
\sigma) V^\T,
\]
where $u_i$, $v_i$ are orthogonal basis in $\mathbb{R}^{p_1}$,
$\mathbb{R}^{p_2}$, respectively, and $\vec\sigma$ is the singular
value vector such that $\sigma_1\geq\sigma_2\geq\cdots\geq\sigma
_p\geq0$. Without loss of generality, we can assume $p \geq kr$,
otherwise we can set the undefined entries of $\sigma$ as $0$.

Consider the singular value vector $\vec\sigma= (\sigma_1, \sigma
_2, \ldots, \sigma_p)$, we note that\break $\vec\sigma_{-\max(kr)}$ satisfies
\begin{eqnarray*}
\llVert \vec\sigma_{-\max(kr)}\rrVert _\infty&\leq& \sigma_{kr},
\\
\llVert \vec\sigma_{-\max(kr)}\rrVert _1 &= & \llVert \vec
\sigma_{-\max(r)}\rrVert _1 - (\sigma_{r+1} +\cdots+
\sigma_{kr})
\\
&\leq& \llVert \vec\sigma_{-\max
(r)}\rrVert _1 - (k-1)r
\sigma_{kr}
\\
&\leq&\llVert \vec\sigma_{\max(r)}\rrVert
_1 - (k-1)r \sigma_{kr}.
\end{eqnarray*}
Denote $\theta= \max \{\sigma_{kr}, (\llVert  \vec\sigma_{\max(r)}\rrVert
_1-r(k-1)\sigma_{kr})/(kr) \}$, by the two inequalities above we
have $\llVert  \vec\sigma_{-\max(kr)}\rrVert  _\infty\leq\theta$ and $\llVert  \vec
\sigma_{-\max(kr)}\rrVert  _1\leq kr\theta$. Now apply Lemma~\ref
{lmsparsemean}, we can get $b^{(i)}\in\mathbb{R}^{p} ,\lambda
_i\geq0 ,i=1,\ldots, N$ such that $\sum_{i=1}^N \lambda_i = 1$,
$\vec\sigma_{-\max(kr)} = \sum_{i=1}^N \lambda_i b^{(i)}$ and
%
\begin{eqnarray}
\label{eqbsparsemean} \operatorname{supp}\bigl(b^{(i)}\bigr) &\subseteq&
\operatorname{supp}(\vec \sigma_{-\max(kr)}),\qquad\bigl\llVert
b^{(i)}\bigr\rrVert _0\leq kr,
\nonumber\\[-8pt]\\[-8pt]\nonumber
\bigl\llVert b^{(i)}\bigr\rrVert _1&=& \llVert \vec
\sigma_{-\max(kr)}\rrVert _1,\qquad\bigl\llVert b^{(i)}
\bigr\rrVert _\infty\leq\theta,
\end{eqnarray}
which leads to
\[
\bigl\llVert b^{(i)}\bigr\rrVert _2 \leq\sqrt{\bigl\llVert
b^{(i)}\bigr\rrVert _1 \cdot\bigl\llVert b^{(i)}
\bigr\rrVert _\infty} \leq \sqrt{\bigl(\llVert \vec\sigma_{\max(r)}
\rrVert _1 - r(k-1)\sigma_{kr}\bigr) \cdot \theta}.
\]
If $\theta= \sigma_{kr}$, we have
\begin{eqnarray*}
\bigl\llVert b^{(i)}\bigr\rrVert _2 & \leq&\sqrt{\bigl(
\llVert \vec\sigma_{\max(r)}\rrVert _1 - r(k-1)
\sigma_{kr}\bigr) \sigma_{kr}}
\\
& \leq&\sqrt{ \biggl(\llVert \vec\sigma_{\max(r)}\rrVert _1 -
r(k-1) \frac{\llVert
\vec\sigma_{\max(r)}\rrVert  _1}{2r(k-1)} \biggr) \frac{\llVert  \vec\sigma
_{\max(r)}\rrVert  _1}{2r(k-1)} }
\\
& \leq&\frac{\llVert  \vec\sigma_{\max(r)}\rrVert  _1}{\sqrt{4r(k-1)}} \leq \frac{\llVert  \vec\sigma_{\max(r)}\rrVert  _2}{\sqrt{4(k-1)}}.
\end{eqnarray*}
If $\theta= (\llVert  \vec\sigma_{\max(r)}\rrVert  _1-r(k-1)\sigma_{kr}) /
(kr)$, we have
\[
\bigl\llVert b^{(i)}\bigr\rrVert _2 \leq\sqrt{
\frac{1}{kr}} \bigl(\llVert \vec\sigma_{\max(r)}\rrVert
_1 - r(k-1) \sigma_{kr}\bigr) \leq\sqrt{
\frac{1}{kr}} \llVert \vec\sigma_{\max
(r)}\rrVert _1 \leq
\frac{\llVert  \vec\sigma_{\max(r)}\rrVert  _2}{\sqrt{k}}.
\]
Since $k\geq2$, we always have $\llVert  b^{(i)}\rrVert  _2\leq\llVert  \vec\sigma_{\max
(r)}\rrVert  _2/\sqrt{k}$. Finally, we define
$B_i = U\operatorname{diag}(b^{(i)}) V^\T$, then the rank of $B_i$ are
all at most $kr$ and $\sum_{i=1}^N \lambda_iB_i = R_{-\max(kr)}$ and
\[
\llVert B_i\rrVert _F = \bigl\llVert b^{(i)}
\bigr\rrVert _2 \leq\llVert \vec\sigma_{\max(r)}\rrVert
_2 / \sqrt {k} = \llVert R_{\max(r)}\rrVert _F/
\sqrt{k}.
\]
Hence,
\begin{eqnarray*}
0 & =& \bigl\llVert \mathcal{X}(R)\bigr\rrVert _1 \geq\bigl\llVert
\mathcal{X}(R_{\max(kr)})\bigr\rrVert _1 - \bigl\llVert
\mathcal{X}(R_{-\max(kr)})\bigr\rrVert _1
\\
& \geq& C_1 \llVert R_{\max(kr)}\rrVert _F - \sum
_{i=1}^N \bigl\llVert \mathcal{X}(\lambda
_i B_i)\bigr\rrVert _1
\\
& \geq& C_1 \llVert R_{\max(r)}\rrVert _F - \sum
_{i=1}^N \lambda_i
C_2 \llVert B_i\rrVert _F
\\
& \geq& C_1 \llVert R_{\max(r)}\rrVert _F -
C_2 \llVert R_{\max(r)}\rrVert _F/\sqrt{k} > 0.
\end{eqnarray*}
Here, we used the RUB condition. The last inequality is due to
$C_2/C_1<\sqrt{k}$ and $R\neq0$ (so $R_{\max(r)}\neq0$). This is a
contradiction, which completes the proof of the theorem.


\subsection{Proof of Theorem \texorpdfstring{\protect\ref{thgaussianl1bound}}{2.2}}
Notice that for $\mathcal{P}$ as standard Gaussian distribution, the
constant $\alpha_{\mathcal{P}}$ [defined as (\ref{eqalphaP})]
equals $1$. We will prove the following more general result than
Theorem \ref{thgaussianl1bound} instead. The proof is provided in
the supplementary material \cite{Supplement}.

\begin{Proposition}\label{prsub-gaussianl1bound}
Suppose $\mathcal{X}\dvtx \mathbb{R}^{p_1\times p_2} \to\mathbb{R}^n$ is
ROP from some variance 1 symmetric sub-Gaussian distribution $\mathcal
{P}$. For integer $k\geq2$, positive $C_1 < \frac{1}{3\alpha
^4_\mathcal{P}}$ [$\alpha_\mathcal{P}$~is defined as (\ref
{eqalphaP})] and $C_2>1$, there exists constants $C$ and $\delta$,
only depending on $\mathcal{P}, C_1, C_2$ but not on $p_1, p_2, r$,
such that if $n\geq Cr(p_1+p_2)$, then with probability at least
$1-e^{-n\delta}$, $\mathcal{X}$ satisfies RUB of order $kr$ and
constants $C_1$ and~$C_2$.
\end{Proposition}

\subsection{Proof of Theorems~\texorpdfstring{\protect\ref{thmixedestimator}}{2.3} 
and  
\texorpdfstring{\protect\ref{thsubgaussiannoise}}{3.1}, 
Proposition \texorpdfstring{\protect\ref{approx-low-rankprop}}{2.1}}
In order to prove the result, we introduce the following technical
lemma as an extension of null space property (Lemma \ref
{lmnullspace}) from exact low-rank into the approximate low-rank setting.

\begin{Lemma}\label{lmnullspaceapproximatelowrank}
Suppose $A_\ast, A\in\mathbb{R}^{p_1\times p_2}$, $R = A_\ast- A$.
If $\llVert  A_\ast\rrVert  _\ast\leq\llVert  A \rrVert  _\ast$, we have
%
\begin{equation}
\label{eqR-maxr} \llVert R_{-\max(r)}\rrVert _\ast\leq\llVert
R_{\max(r)}\rrVert _\ast+ 2 \llVert A_{-\max(r)}\rrVert
_\ast.
\end{equation}
\end{Lemma}

The following two lemmas described the separate effect of constraint
$\mathcal{Z}_1 = \{z\dvtx  \llVert  z\rrVert  _1/n\leq\lambda_1\}$ and $\mathcal{Z}_2
= \{z\dvtx  \llVert  \mathcal{X}^\ast(z)\rrVert  \leq\lambda_2 \}$ on the estimator.

\begin{Lemma}\label{lmnoisyl1}
Suppose $\mathcal{X}$ satisfies RUB condition of order $kr$ with
constants $C_1, C_2$ such that $C_1>C_2/\sqrt{k}$. Assume that $A_\ast
, A\in\mathbb{R}^{p_1\times p_2}$ satisfy $\llVert  A_\ast\rrVert  _\ast\leq\llVert  A \rrVert
_\ast$, $\llVert  \mathcal{X}(A_\ast- A)\rrVert  _1/n \leq\lambda_1$. Then we have
\[
\llVert A_\ast- A\rrVert _F \leq\frac{2}{C_1 - C_2/\sqrt{k}}
\lambda_1 + \biggl(\frac{3}{\sqrt{k}C_1/C_2-1} + \frac{1}{\sqrt{k-1}} \biggr)
\frac{\llVert
A_{-\max(r)}\rrVert  _\ast}{\sqrt{r}}.
\]
\end{Lemma}

%
\begin{Lemma}\label{lmnoisyDS}
Suppose $\mathcal{X}$ satisfies RUB condition of order $kr$ with
constants $C_1, C_2$ such that $C_1>C_2/\sqrt{k}$. Assume that $ \hat A{}^{\mathrm{DS}}$ satisfies $\llVert  \mathcal{X}^\ast\mathcal{X}(A_\ast- A)\rrVert   \leq
\lambda_2$. Then we have
\begin{eqnarray*}
\llVert A_\ast-A\rrVert _F &\leq& \frac{4}{(C_1 - C_2/\sqrt{k})^2}\cdot
\frac{\sqrt{r}\lambda_2}{n}
\\
&&{} + \biggl(\frac{5}{\sqrt{k}C_1/C_2 - 1} + \frac{1}{\sqrt{k-1}} + 1 \biggr)
\frac{\llVert  A_{-\max(r)}\rrVert  _\ast
}{\sqrt{r}}.
\end{eqnarray*}
\end{Lemma}

The proof of Lemmas \ref{lmnullspaceapproximatelowrank}, \ref
{lmnoisyl1} and \ref{lmnoisyDS} are listed in the supplementary material \cite{Supplement}. Now
we prove Theorem \ref{thmixedestimator} and Proposition \ref
{approx-low-rankprop}. We only need to prove Proposition \ref
{approx-low-rankprop} since Theorem \ref{thmixedestimator} is a
special case of Proposition \ref{approx-low-rankprop}. By Lemmas \ref{lml1norm}~and~\ref{lmDSbound}, we have
\begin{eqnarray*}
&& P_z\bigl(\llVert z\rrVert _1\leq\sigma n \bigr) \leq
\frac{9}{n},
\\
&& P_{\mathcal{X}, z} \bigl(\bigl\llVert \mathcal{X}^\ast(z)\bigr\rrVert
\geq\sigma \bigl(12(p_1 + p_2)\sqrt{\log n} + 6
\sqrt{2(p_1+p_2)n} \bigr) \bigr)
\\
&&\qquad \leq P_{\mathcal{X}} \bigl(\bigl\llVert \mathcal{X}^\ast(z)\bigr
\rrVert \geq \bigl(6(p_1+p_2)\llVert z\rrVert
_\infty+ 6\sqrt{p_1+p_2}\llVert z\rrVert
_2 \bigr) \bigr)
\\
&&\quad\qquad{} + P_{z} \bigl(\llVert z\rrVert _\infty\geq2\sigma\sqrt{
\log n} \bigr) + P_{z} \bigl(\llVert z\rrVert _2 \geq
\sigma\sqrt{2n} \bigr)
\\
&&\qquad \leq 2\exp\bigl(-(2-\log7) (p_1 + p_2)\bigr) +
\frac{1}{n\sqrt{2\pi\log n}} + \frac{1}{n}.
\end{eqnarray*}
Here, $P_{\mathcal{X}}$ ($P_z$ or $P_{\mathcal{X}, z}$) means the
probability with respect to $\mathcal{X}$ [$z$ or $(\mathcal{X},
z)$]. Hence, we have
\[
P(z\in\mathcal{Z}_1\cap\mathcal{Z}_2) \geq1 - 2\exp
\bigl(-(2-\log 7) (p_1 + p_2)\bigr) - \frac{11}{n}.
\]
Under the event that $z\in\mathcal{Z}_1\cap\mathcal{Z}_2$, $A$ is
in the feasible set of the programming (\ref
{eqnuclearnormminimization2}), which implies $\llVert  \hat A\rrVert  _\ast\leq\Vert A\Vert_\ast$ by the definition of $\hat A$. Moreover, we have
\begin{eqnarray*}
\bigl\llVert \mathcal{X}(\hat A - A)\bigr\rrVert _1/n &\leq& \bigl
\llVert y - \mathcal{X}(A)\bigr\rrVert _1/n + \bigl\llVert y -
\mathcal{X}(\hat A)\bigr\rrVert _1/n
\\
&\leq& \llVert z\rrVert _1/n + \bigl\llVert y - \mathcal{X}(\hat A)
\bigr\rrVert _1/n\leq2 \sigma,
\\
\bigl\llVert \mathcal{X}^\ast\mathcal{X}(\hat A- A)\bigr\rrVert &\leq&
\bigl\llVert \mathcal {X}^\ast\bigl(y - \mathcal{X}(\hat A)\bigr)\bigr
\rrVert + \bigl\llVert \mathcal{X}^\ast\bigl(y - \mathcal{X}(A)\bigr)
\bigr\rrVert
\\
&\leq& \bigl\llVert \mathcal{X}^\ast\bigl(y - \mathcal{X}(\hat A)\bigr)
\bigr\rrVert + \bigl\llVert \mathcal {X}^\ast(z)\bigr\rrVert \leq2\eta.
\end{eqnarray*}
On the other hand, suppose $k = 10$, by Theorem \ref
{thgaussianl1bound}, we can have find a uniform constant $C$ and
$\delta$ such that if $n\geq Crk(p_1+p_2)$, $\mathcal{X}$ satisfies
RUB of order $10r$ and constants $C_1 = 0.32, C_2 = 1.02$ with
probability at least $1-e^{-n\delta'}$. Hence, we have $D (= Ck)$ and
$\delta'$ such that if $n\geq Dr(p_1+p_2)$, $\mathcal{X}$ satisfies
RUB of order $10r$ and constants $C_1, C_2$ satisfying $C_2/C_1 < \sqrt
{10}$ with probability at least $1-e^{-n\delta'}$.

Now under the event that:
\begin{longlist}[2.]
\item[1.] $\mathcal{X}$ satisfies RUB of order $10r$ and constants $C_1,
C_2$ satisfying $C_2/C_1<\sqrt{10}$,
\item[2.] $z\in\mathcal{Z}_1\cap\mathcal{Z}_2$,
\end{longlist}
apply Lemmas \ref{lmnoisyl1}~and~\ref{lmnoisyDS} with
$A_\ast= \hat A$, we can get (\ref{eqapproxlowrankinequality}).
The probability that these two events both happen is at least\vspace*{1pt}
$1 - 2\exp(-(2-\log7)(p_1 + p_2)) - \frac{11}{n} - \exp(-\delta' n)$.
Set $\delta= \min(2-\log7, \delta')$, we finished the proof of
Proposition \ref{approx-low-rankprop}.

For Theorem \ref{thsubgaussiannoise}, the proof is similar. We
apply the latter part of Lemmas~\ref{lml1norm}~and~\ref{lmDSbound} and get
\begin{eqnarray*}
&& P(z\notin\mathcal{Z}_1 \cap\mathcal{Z}_2)
\\
&&\qquad \leq P\bigl(\llVert z\rrVert _1/n > 6\tau\bigr)
\\
&&\quad\qquad{}  + P \bigl(\bigl
\llVert \mathcal{X}(z)\bigr\rrVert > \tau \alpha^2_{\mathcal{P}}
\bigl(6\sqrt{6n(p_1+p_2)} + 12\sqrt{\log
n}(p_1 + p_2) \bigr) \bigr)
\\
&&\qquad \leq P\bigl(\llVert z\rrVert /n > 6 \tau\bigr) + P\bigl(\llVert z\rrVert
_2 > \sqrt{6n}\tau\bigr) + P\bigl(\llVert z\rrVert _\infty>
2\sqrt{\log n} \tau\bigr)
\\
&&\quad\qquad{}  + P_{\mathcal{X}}\bigl(\bigl\llVert \mathcal{X}(z)\bigr\rrVert >
\alpha_{\mathcal
{P}}^2\bigl(6(p_1+p_2)
\llVert z\rrVert _\infty+ 6\sqrt{p_1+p_2}\llVert
z\rrVert _2 \bigr)\bigr)
\\
&&\qquad \leq \exp \bigl(-n(6 - 2\sqrt{2\pi})^2/2 \bigr) + \exp(-n/12)
\\
&&\quad\qquad{} +\frac{2}{n} + 2\exp\bigl(-(2-\log7) (p_1 + p_2)
\bigr).
\end{eqnarray*}
Besides, we choose $k > (3\alpha^4_{\mathcal{P}})^2$, then we can
find $C_1 < 1/(3\alpha_{\mathcal{P}}^4)$ and $C_2 > 1$ such that $C_2
/ C_1 < \sqrt{k}$. Apply Proposition \ref{prsub-gaussianl1bound},
there exists $C,\delta'$ only depending on $\mathcal{P}$, $C_1, C_2$
such that if $n\geq Ckr(p_1 + p_2)$, $\mathcal{X}$ satisfies RUB of
order $kr$ with constants $C_1$ and $C_2$ with probability at least $1
- \exp(-\delta'(p_1 + p_2))$. Note that $C_1, C_2$ only depends on
$\mathcal{P}$, we can conclude that there exist constants $D (= Ck),
\delta'$ only depending on $\mathcal{P}$ such that if $n\geq Dr(p_1 +
p_2)$, $\mathcal{X}$ satisfies RUB of order $kr$ with constants $C_1,
C_2$ satisfying $C_2/C_1 \leq\sqrt{k}$.

Similarly, to the proof of Proposition \ref{approx-low-rankprop},
under the event that:
\begin{longlist}[2.]
\item[1.] $\mathcal{X}$ satisfies RUB of order $kr$ and constants $C_1,
C_2$ satisfying $C_2/C_1<\sqrt{k}$,
\item[2.] $z\in\mathcal{Z}_1\cap\mathcal{Z}_2$,
\end{longlist}
we can get (\ref{eqsubgaussiannoiseinequality}) (we shall note
that $W_1$ depends on $\mathcal{P}$, so its value can also depend on
$\alpha_{\mathcal{P}}$). The probability that those events happen is
at least $1 - 2/n - 5\exp(-\delta(p_1+p_2))$ for $\delta\leq\min
((6-2\sqrt{2\pi})^2/2, 1/12, 2-\log7, \delta')$.

\subsection{Proof of Theorem \texorpdfstring{\protect\ref{thlowerbound}}{2.4}}
Without loss of generality, we assume that \mbox{$p_1 \leq p_2$}. We consider
the class of rank-$r$ matrices
\[
\mathcal{F}_c = \bigl\{A\in\mathbb{R}^{p_1\times p_2}\dvtx
A_{ij} = 0, \mbox { whenever } i\geq r+1\bigr\}
\]
namely the matrices with all nonzero entries in the first $r$ rows. The
model (\ref{eqmodel}) become
\[
y_i = \beta_{1\dvtx r}^{(i)T}A_r
\gamma^{(i)}+z_i,\qquad i = 1,\ldots, n,
\]
where $\beta^{(i)}_{1\dvtx r}$ is the vector of the first to the $r$th
entries of $\beta^{(i)}$. Note that this is a linear regression model
with variable $A_r \in\mathbb{R}^{r\times p_2}$, by Lemma 3.11 in
\cite{CandesOracle}, we have
%
\begin{eqnarray}
\inf_{\hat A}\sup_{A\in\mathcal{F}_c} E\bigl\llVert \hat
A(y) - A\bigr\rrVert _F^2 &=& \sigma^2 \operatorname{trace} \bigl[\bigl(\mathcal{X}_r^\ast\mathcal
{X}_r\bigr)^{-1} \bigr],
\\
\inf_{\hat A}\sup_{A\in\mathcal{F}_c} E\bigl\llVert \hat
A(y) - A\bigr\rrVert _F^2 &=& \infty\qquad\mbox{when $
\mathcal{X}_r^\ast\mathcal{X}_r$ is singular,}
\end{eqnarray}
where $\mathcal{X}_r\dvtx \mathbb{R}^{r\times p_2}\to\mathbb{R}^n$ is
the $\mathcal{X}$ constrained on $\mathcal{F}_c$, Then $\mathcal
{X}_r$ sends $A_r$ to $ (\beta^{(1)}_{1\dvtx r}A_r\gamma^{(1)},
\ldots, \beta^{(n)}_{1\dvtx r}A_r\gamma^{(n)} )^\T$. When $n <
p_2r$, $\mathcal{X}_r$ is singular, hence we have~(\ref{eqlowerboundinfty}).

When $n\geq p_2r$, we can see in order to show (\ref
{eqlowerboundE}), we only need to show $\operatorname
{trace}(\mathcal{X}_r^\ast\mathcal{X}_r)\geq\frac{p_2 r}{2n} $
with probability at least $1-26n^{-1}$. Suppose the singular value of
$\mathcal{X}_r$ are $\sigma_i(\mathcal{X}_r)$, $i=1,\ldots, rp_2$, then
$\operatorname{trace}(\mathcal{X}_r^\ast\mathcal{X}_r) = \sum_{i=1}^{p_2r}\sigma^{-2}(\mathcal{X}_r)$.

Suppose $\mathcal{X}$ is ROP while $B \in\mathbb{R}^{r\times p_2}$
is i.i.d. standard Gaussian random matrix (both $\mathcal{X}$ and
$B_r$ are random). Then by some calculation, we can see
\[
E_{B, \mathcal{X}_r}\bigl\llVert \mathcal{X}_r(B)\bigr\rrVert
_2^2 = n E_{B, \beta, \gamma
} \bigl(\beta_{1\dvtx r}^\T B
\gamma \bigr)^2 = n\sum_{j=1}^r
\sum_{k=1}^{p_2}E (\beta_j
B_{jk}\gamma_k )^2 = np_2r.
\]
Note (0.20) 
in the proof of Lemma \ref{lmmgf} in the supplementary material \cite{Supplement}, we know $E
(\beta_{1\dvtx r}^{(i)T}B\gamma^{(i)}\Vert _2^4\mid B ) \leq9\llVert  B\rrVert
_F^4$. Hence,
\begin{eqnarray*}
E\bigl\llVert \mathcal{X}_r(B)\bigr\rrVert _2^4
&=&  \sum_{i=1}^n E \bigl(
\beta_{1\dvtx r}^{(i)T}B\gamma^{(i)} \bigr)^4
\\
&&{} +
2\sum_{1\leq i<l\leq n}E\sum_{j=1}^n
\bigl(\beta_{1\dvtx r}^{(i)T}B\gamma ^{(i)}
\bigr)^2\cdot E\sum_{j=1}^n
\bigl(\beta_{1\dvtx r}^{(l)T}B\gamma ^{(l)}
\bigr)^2
\\
&=&  n \cdot9E\llVert B\rrVert _F^4 + n(n-1)
(p_2r)^2
\\
&=& 9n E \bigl(\chi ^2(p_2r)
\bigr)^2 + n(n-1)p_2^2r^2
\\
&=&  9n\bigl(p_2^2r^2 + 2p_2r
\bigr) + n(n-1)p_2^2r^2
\\
&=& n^2p_2^2r^2 +
2np_2r(4p_2r+9) \leq n^2p_2^2r^2
+ 26np_2^2r_2^2.
\end{eqnarray*}
Besides,
\begin{eqnarray*}
E\bigl\llVert \mathcal{X}_r(B_r)\bigr\rrVert
_2^2 &=& E \bigl(E \bigl(\bigl\llVert \mathcal{X}_r(B_r)
\bigr\rrVert _2^2\mid \mathcal{X}_r \bigr)
\bigr) =E \Biggl(\sum_{i=1}^{rp_2}
\sigma_i^2(\mathcal{X}_r) \Biggr),
\\
E\bigl\llVert \mathcal{X}_r(B_r)\bigr\rrVert
_2^4 &=& E \bigl(E \bigl(\bigl\llVert \mathcal{X}_r(B_r)
\bigr\rrVert _2^4\mid\mathcal{X}_r \bigr)
\bigr)
\\
&= &E \Biggl(\sum_{i=1}^{rp_2}3
\sigma_i^4(\mathcal{X}_r) + 2\sum
_{1\leq i<j\leq rp_2}\sigma_i^2(
\mathcal{X}_r)\sigma_j^2(\mathcal
{X}_r) \Biggr)
\\
&\geq& E \Biggl(\sum_{i=1}^{rp_2}
\sigma_i ^2(\mathcal {X}_r)^2
\Biggr)^2.
\end{eqnarray*}
Hence,
\begin{eqnarray*}
E \Biggl(\sum_{i=1}^{rp_2}
\sigma_i ^2(\mathcal{X}_r)^2
\Biggr) &=& np_2r,
\\
\operatorname{Var} \Biggl(\sum_{i=1}^{rp_2}
\sigma_i ^2(\mathcal{X}_r)^2
\Biggr) &=& E \Biggl(\sum_{i=1}^{rp_2}
\sigma_i ^2(\mathcal{X}_r)^2
\Biggr)^2 - \Biggl(E\sum_{i=1}^{rp_2}
\sigma_i ^2(\mathcal{X}_r)^2
\Biggr)^2 \leq26np_2^2r^2.
\end{eqnarray*}
Then by Chebyshev's inequality, we have
%
\begin{equation}
\label{eqineqsigmai^2X} \sum_{i=1}^{rp_2}
\sigma_i^2 (\mathcal{X}_r)
\leq2np_2r
\end{equation}
with probability at least $1- \frac{26np_2^2r^2}{(npr)^2} = 1 - \frac{26}{n}$.
By Cauchy--Schwarz's inequality, we have
\[
\operatorname{trace} \bigl(\bigl(\mathcal{X}_r^\ast
\mathcal {X}_r\bigr)^{-1} \bigr) = \sum
_{i=1}^{rp_2} \sigma_i^{-2}(
\mathcal{X}_r) \geq\frac{(p_2r)^2}{\sum_{i=1}^{rp_2} \sigma_i^{2}(\mathcal{X}_r)}.
\]
Therefore, we have
\[
\operatorname{trace} \bigl(\bigl(\mathcal{X}_r^\ast
\mathcal {X}_r\bigr)^{-1} \bigr) \geq\frac{p_2r}{2n}
\]
with probability at least $1- 26/n$, which shows (\ref{eqlowerboundE}).

Finally, we consider (\ref{eqlowerboundP}). Suppose inequality
(\ref{eqineqsigmai^2X}) holds, then
%
\begin{eqnarray}\label{eqlowerboundnumber}
&& \bigl\llvert \bigl\{i\dvtx  \sigma_i^{2}(X_r)
\geq4n\bigr\}\bigr\rrvert \leq\frac{p_2r}{2}\nonumber
\\
&&\qquad\Rightarrow\quad \biggl\llvert \biggl\{i\dvtx  \sigma_i^{-2}(X_r)
\leq\frac{1}{4n}\biggr\}\biggr\rrvert \geq\frac{p_2r}{2}
\\
&&\qquad\Rightarrow\quad \biggl\llvert \biggl\{i\dvtx  \sigma_i^{-2}(X_r)
\geq\frac{1}{4n}\biggr\}\biggr\rrvert \geq\frac{p_2r}{2}.\nonumber
\end{eqnarray}
By Lemma 3.12 in \cite{CandesOracle}, we know
\begin{eqnarray*}
&&  \inf_{\hat A}\sup_{A\in\mathcal{F}_c} P_z
\biggl(\llVert \hat A - A\rrVert _F^2 \geq
\frac{p_2r\sigma^2}{16n} \biggr)
\\
&&\qquad =  \inf_{\hat A}\sup_{A\in\mathcal{F}_c} E_z
1_{\{x \geq p_2r\sigma^2/16n\}}\bigl(\llVert \hat A - A\rrVert _F^2
\bigr)
\\
&&\qquad =  E_z 1_{\{x\geq p_2r\sigma^2/16n\} }\bigl(\bigl\llVert \bigl(\mathcal
{X}_r^\ast\mathcal{X}_r\bigr)^{-1}
\mathcal{X}_r^\ast(z)\bigr\rrVert _F^2
\bigr)
\\
&&\qquad =  P_z \biggl(\bigl\llVert \bigl(\mathcal{X}_r^\ast
\mathcal{X}_r\bigr)^{-1}\mathcal {X}_r^\ast(z)
\bigr\rrVert _F^2\geq\frac{p_2r\sigma^2}{16n} \biggr),
\end{eqnarray*}
where $1_{\{x\geq p_2r\sigma^2/16n\}}(\cdot)$ is the
indicator function. Note that when $z\stackrel{\mathrm{i.i.d.}}{\sim}\break  N(0, \sigma
^2)$, $\llVert  (\mathcal{X}_r^\ast\mathcal{X}_r)^{-1}\mathcal{X}_r^\ast
(z)\rrVert  _F^2$ is identical distributed as
$\sum_{i=1}^{rp_2}\frac{y_i^2}{\sigma_i^2(\mathcal{X}_r)}$,
where $y_1,\ldots,\break  y_{rp_2}\stackrel{\mathrm{i.i.d.}}{\sim} N(0, \sigma^2)$, hence,
\begin{eqnarray*}
&& P \biggl(\bigl\llVert \bigl(\mathcal{X}_r^\ast
\mathcal{X}_r\bigr)^{-1}\mathcal {X}_r^\ast(z)
\bigr\rrVert _F^2 \leq\frac{p_2r\sigma^2}{16n} \biggr)
\\
&&\qquad =  P \Biggl(\sum_{i=1}^{rp_2}
\frac{y_i^2}{\sigma_i^2(\mathcal
{X}_r)} \leq\frac{p_2r\sigma^2}{16n} \Biggr)
\\
&&\qquad \leq P \biggl(\sum_{i\dvtx \sigma_i^{-2}(\mathcal{X}_r)\geq1/(4n)} y_i^2
\sigma_i^{-2}(\mathcal{X}_r) \leq
\frac{p_2r\sigma
^2}{16n} \biggr)
\\
&&\qquad \leq P \biggl(\sum_{i\dvtx  \sigma_i^{-2}(\mathcal{X}_r)\geq1/(4n)} \frac{y_i^2}{4n} \leq
\frac{p_2r\sigma^2}{16n} \biggr) \leq P \biggl(\chi^2\biggl(\biggl\lceil
\frac{rp_2}{2}\biggr\rceil\biggr) \leq\frac{p_2r}{4} \biggr)
\\
&&\qquad \leq \exp \biggl(-\frac{rp_2}{32} \biggr).
\end{eqnarray*}
The last inequality is due to the tail bound of $\chi^2$ distribution
given by Lemma 1 in~\cite{Laurent}; the second last inequality is due
to (\ref{eqlowerboundnumber}). In summary, when (\ref
{eqineqsigmai^2X}) holds, we have
\[
\inf_{\hat A}\sup_{A\in\mathcal{F}_c} P_z
\biggl(\llVert \hat A - A\rrVert _F^2 \geq
\frac{p_2r\sigma^2}{16n} \biggr) \leq\exp\biggl( - \frac{rp_2}{32}\biggr).
\]
Finally,\vspace*{1pt} since $p_2 \geq(p_1+p_2)/2$, we showed that with probability
at least $1-26n^{-1}$, $\mathcal{X}$ satisfies (\ref
{eqlowerboundP}). 

\subsection{Proof of Theorem \texorpdfstring{\protect\ref{thcovariance}}{4.1}}
 
We first introduce the following lemma about the upper bound of $\llVert  z\rrVert
_1, \llVert  z\rrVert  _2, \llVert  z\rrVert  _\infty$.

\begin{Lemma}\label{lmcovarianceerrorbound}
Suppose $z$ is defined as (\ref{eqyandz}), then for constants
$C_1>\sqrt{2}$, $M_1>1$, we have
%
\begin{eqnarray}
\label{eqcovarianceinequalityl1} P \Biggl(\llVert z\rrVert _1/n \leq
\frac{C_1}{n}\sum_{i=1}^n
\xi_{i}^2 \Biggr)&\geq& 1-\frac{9C_1^2+6}{n(C_1-\sqrt{2})^2},
\nonumber\\[-8pt]\\[-8pt]\nonumber
P \Biggl(\frac{C_1}{n}\sum_{i=1}^n
\xi_{i}^2 \leq M_1C_1\llVert
\Sigma\rrVert _\ast \Biggr)&\geq& 1 - \frac{9}{n(M_1-1)^2};
\end{eqnarray}
%
for constants $C_2> 1$, $M_2>9$,
%
\begin{eqnarray}
\label{eqcovarianceinequalityl2}  P \biggl(\llVert z\rrVert _2^2/n \leq
\frac{C_2^2\sum_{i=1}^n\xi
_{i}^4}{n} \biggr)&\geq& 1-\frac{105 (105 C_2^4 + 60
)}{n(3C_2^2 - 2)^2},
\nonumber\\[-8pt]\\[-8pt]\nonumber
P \biggl(\frac{C_2^2\sum_{i=1}^n\xi_{i}^4}{n} \leq M_2C_2^2
\llVert \Sigma \rrVert _\ast^2 \biggr)&\geq& 1 -
\frac{105^2}{n(M_1 - 9)^2};
\end{eqnarray}
%
for constants $C_3>1$, $M_3>1$,
%
\begin{eqnarray}\label{eqcovarianceinequalityinfty}
\qquad &&  P \Bigl(\llVert z\rrVert _\infty\leq C_3
\log n \max_{1\leq i\leq n}\xi _{i}^2 \Bigr) \geq 1
- \frac{2}{\sqrt{2\pi C_3\log n}},\nonumber
\\
&& P \Bigl(C_3\log n \max_{1\leq i\leq n}
\xi_{i}^2\leq2C_3 M_3 \log
^2 n \bigl(\sqrt{\llVert \Sigma\rrVert _\ast} +
\sqrt{2M_3 \log n\llVert \Sigma\rrVert } \bigr)^2 \Bigr)
\\
&&\qquad \geq 1-2n^{-M_3+1}.\nonumber
\end{eqnarray}
%
\end{Lemma}

The proof of Lemma \ref{lmcovarianceerrorbound} is listed in the
supplementary material \cite{Supplement}. The rest of the proof is basically the same as Proposition
\ref{prsymmetricmixedestimator}. Suppose $\mathcal{X}_1, \mathcal
{X}_2$ and $\tilde z$ are given by (0.36), (0.37) and (0.39) in the
supplementary material \cite{Supplement}, then $\mathcal{X}_1$, $\mathcal{X}_2$ are ROP. By Lemma \ref{lmDSbound},
%
\begin{eqnarray}\label{eqcovariancenormX1}
\bigl\llVert \mathcal{X}_1^\ast(\tilde z)
\bigr\rrVert &\leq& 6 \bigl(2p\llVert \tilde z\rrVert _\infty+ \sqrt{2p}
\llVert \tilde z\rrVert _2 \bigr),
\\
\label{eqcovariancenormX2} \bigl\llVert \mathcal{X}_2^\ast(\tilde z)
\bigr\rrVert &\leq& 6 \bigl(2p\llVert \tilde z\rrVert _\infty+ \sqrt{2p}
\llVert \tilde z\rrVert _2 \bigr)
\end{eqnarray}
with probability at least $1 - 4\exp (-2(2 - \log7)p )$.
Hence, there exists $\delta>0$ such that
\begin{eqnarray*}
&& P \bigl(\Sigma_0 \mbox{ is NOT in the feasible set of
(\ref
{eqminimizationcovariance})} \bigr)
\\
&&\qquad  = P \bigl(\llVert z\rrVert _1/n > \eta
_1\mbox{ or }\bigl\llVert \tilde{\mathcal{X}}^\ast(
\tilde{z})\bigr\rrVert > \eta_2 \bigr)
\\
&&\qquad \leq P \Biggl(\llVert z\rrVert _1/n > \frac{c_1}{n}\sum
_{i=1}^n \xi_i^2
\Biggr) + P \Bigl(\llVert \tilde{z}\rrVert _\infty> 2c_3\log
n \max_{1\leq i\leq n}\xi _{i}^2 \Bigr)
\\
&&\quad\qquad{}  + P \Biggl(\llVert \tilde z\rrVert _2 > c_2\sqrt{2
\sum_{i=1}^n\xi _i^4}
\Biggr)
\\
&&\quad\qquad{} + P \bigl(\bigl\llVert \tilde{\mathcal{X}}^\ast(z)\bigr\rrVert > 24p
\llVert \tilde z\rrVert _\infty + 12\sqrt{2p}\llVert \tilde z\rrVert
_2 \bigr)
\\
&&\qquad \leq P \Biggl(\llVert z\rrVert _1/n > \frac{c_1}{n}\sum
_{i=1}^n \xi_i^2
\Biggr) + P \Bigl(\llVert z\rrVert _\infty> c_3\log n \max
_{1\leq i\leq n} \xi_{i}^2 \Bigr)
\\
&&\quad\qquad{} + P \Biggl(
\llVert z\rrVert _2 > c_2\sqrt{\sum
_{i=1}^n\xi_i^4} \Biggr)
\\
&&\quad\qquad{} + P \bigl(\bigl\llVert \mathcal{X}^\ast_1 (\tilde{z})
\bigr\rrVert > 12p\llVert \tilde z\rrVert _\infty+ 6\sqrt{2p}\llVert
\tilde z\rrVert _2 \bigr)
\\
&&\quad\qquad{}  + P \bigl(\bigl\llVert \mathcal
{X}^\ast_2 (z)\bigr\rrVert > 12p\llVert \tilde z\rrVert
_\infty+ 6\sqrt{2p}\llVert \tilde z\rrVert _2 \bigr)
\\
&&\qquad \leq O(1/n) + 4\exp \bigl(-2(2 - \log7)p \bigr) + \frac{2}{\sqrt
{2\pi c_3\log n}}.
\end{eqnarray*}
Here, we used the fact that $\tilde{\mathcal{X}}^\ast= \mathcal
{X}_1^\ast+\mathcal{X}_2^\ast$,
\begin{eqnarray*}
\llVert \tilde{z}\rrVert _2 &=& \sqrt{\sum
_{i=1}^{\lfloor n/2\rfloor}(z_{2i-1}- z_{2i})^2}
\leq\sqrt{\sum_{i=1}^{\lfloor n/2\rfloor} 2
\bigl(z_{2i-1}^2 + z_{2i}^2\bigr)}
\leq\sqrt{2}\llVert z\rrVert _2,
\\
\llVert \tilde{z}\rrVert _\infty &=& \max_{i} \llvert
z_{2i-1} - z_{2i}\rrvert \leq2\max_i
\llvert z_i\rrvert \leq2\llVert z\rrVert _\infty.
\end{eqnarray*}
Similarly to the proof of Proposition \ref
{prsymmetricmixedestimator}, since $\mathcal{X}_1$ is ROP, there
exists constants $D$ and $\delta'$ such that if $n\geq Drp$, $\mathcal
{X}_1$ satisfies RUB of order $10k$ with constants $C_1, C_2$
satisfying $C_2/C_1<\sqrt{10}$ with probability at least $1 -
e^{-n\delta'}$.

Now under the event that:
\begin{longlist}[3.]
\item[1.] $A$ is feasible in (\ref{eqminimizationcovariance}),
\item[2.] $\mathcal{X}_1$ satisfies RUB of order $10k$ with constants
$C_1, C_2$ satisfying $C_2/C_1<\sqrt{10}$,
\item[3.]  the latter part of (\ref{eqcovarianceinequalityl1}), (\ref
{eqcovarianceinequalityl2}) and (\ref
{eqcovarianceinequalityinfty}) hold for some $M_1>1$, $M_2>9$, $M_3>2$,
\end{longlist}
we can prove (\ref{eqcovarianceineq}) similarly as the proof of
Proposition \ref{prsymmetricmixedestimator}, which we omit the
proof here. 

\section*{Acknowledgments}
We thank the Associate Editor and the referees for their thorough and
useful comments which have helped to improve the presentation of the paper.

\begin{supplement}
\sname{Supplement to ``ROP: Matrix recovery via rank-one projections''} 
\slink[doi]{10.1214/14-AOS1267SUPP} 
\sdatatype{.pdf}
\sfilename{AOS1267\_supp.pdf}
\sdescription{We prove the technical lemmas used in the proofs of the main results in
this supplement. The proofs rely on results in \mbox{\cite{Laurent,Cail1,Vershynin,RechtMatrix,CandesOracle,WangLi}} and \cite{Oymak}.}
\end{supplement}


%

\printaddresses

\begin{thebibliography}{43}

\bibitem{Alquier}
%
\begin{barticle}[auto:parserefs-M02]
\bauthor{\bsnm{Alquier},~\bfnm{P.}\binits{P.}},
\bauthor{\bsnm{Butucea},~\bfnm{C.}\binits{C.}},
\bauthor{\bsnm{Hebiri},~\bfnm{M.}\binits{M.}} \AND
\bauthor{\bsnm{Meziani},~\bfnm{K.}\binits{K.}}
(\byear{2013}).
\btitle{Rank penalized estimation of a quantum system}.
\bjournal{Phys. Rev. A}.
\bvolume{88}
\bpages{032133}.
\end{barticle}
%
\bptok{imsref}%
\endbibitem

\bibitem{Andrews}
%
\begin{barticle}[auto:parserefs-M02]
\bauthor{\bsnm{Andrews},~\bfnm{H.~C.}\binits{H.~C.}} \AND
\bauthor{\bsnm{Patterson},~\bfnm{C.~L.}\binits{C.~L.} \bsuffix{III}}
(\byear{1976}).
\btitle{Singular value decomposition (SVD) image coding}.
\bjournal{IEEE Trans. Commun.}
\bvolume{24}
\bpages{425--432}.
\end{barticle}
%
\bptok{imsref}%
\endbibitem

\bibitem{Basri}
%
\begin{barticle}[auto:parserefs-M02]
\bauthor{\bsnm{Basri},~\bfnm{R.}\binits{R.}} \AND
\bauthor{\bsnm{Jacobs},~\bfnm{D.~W.}\binits{D.~W.}}
(\byear{2003}).
\btitle{Lambertian reflectance and linear sub-spaces}.
\bjournal{IEEE Trans. Pattern Anal. Mach. Intell.}
\bvolume{25}
\bpages{218--233}.
\end{barticle}
%
\bptok{imsref}%
\endbibitem

\bibitem{Birnbaum}
%
\begin{barticle}[mr]
\bauthor{\bsnm{Birnbaum},~\bfnm{Aharon}\binits{A.}},
\bauthor{\bsnm{Johnstone},~\bfnm{Iain~M.}\binits{I.~M.}},
\bauthor{\bsnm{Nadler},~\bfnm{Boaz}\binits{B.}} \AND
\bauthor{\bsnm{Paul},~\bfnm{Debashis}\binits{D.}}
(\byear{2013}).
\btitle{Minimax bounds for sparse PCA with noisy high-dimensional data}.
\bjournal{Ann. Statist.}
\bvolume{41}
\bpages{1055--1084}.
\bid{doi={10.1214/12-AOS1014}, issn={0090-5364}, mr={3113803}}
\end{barticle}
%
\bptok{imsref}%
\endbibitem



\bibitem{CaiMaWu1}
%
\begin{barticle}[mr]
\bauthor{\bsnm{Cai},~\bfnm{T.~Tony}\binits{T.~T.}},
\bauthor{\bsnm{Ma},~\bfnm{Zongming}\binits{Z.}} \AND
\bauthor{\bsnm{Wu},~\bfnm{Yihong}\binits{Y.}}
(\byear{2013}).
\btitle{Sparse PCA: Optimal rates and adaptive estimation}.
\bjournal{Ann. Statist.}
\bvolume{41}
\bpages{3074--3110}.
\bid{doi={10.1214/13-AOS1178}, issn={0090-5364}, mr={3161458}}
\end{barticle}
\bptok{imsref}%
\endbibitem

\bibitem{CaiMaWu2}
%
\begin{barticle}[auto]
\bauthor{\bsnm{Cai},~\bfnm{T.~Tony}\binits{T.~T.}},
\bauthor{\bsnm{Ma},~\bfnm{Zongming}\binits{Z.}} \AND
\bauthor{\bsnm{Wu},~\bfnm{Yihong}\binits{Y.}}
(\byear{2014}).
\btitle{Optimal estimation and rank detection for sparse spiked covariance matrices}.
\bjournal{Probab. Theory Related Fields}.
\bnote{To appear}.
\end{barticle}
%
\bptok{imsref}%
\endbibitem

\bibitem{Cail1}
%
\begin{barticle}[mr]
\bauthor{\bsnm{Cai},~\bfnm{T.~Tony}\binits{T.~T.}},
\bauthor{\bsnm{Xu},~\bfnm{Guangwu}\binits{G.}} \AND
\bauthor{\bsnm{Zhang},~\bfnm{Jun}\binits{J.}}
(\byear{2009}).
\btitle{On recovery of sparse signals via {$\ell\sb1$} minimization}.
\bjournal{IEEE Trans. Inform. Theory}
\bvolume{55}
\bpages{3388--3397}.
\bid{doi={10.1109/TIT.2009.2021377}, issn={0018-9448}, mr={2598028}}
\end{barticle}
%
\bptok{imsref}%
\endbibitem

\bibitem{CaiZhang}
%
\begin{barticle}[mr]
\bauthor{\bsnm{Cai},~\bfnm{T.~Tony}\binits{T.~T.}} \AND
\bauthor{\bsnm{Zhang},~\bfnm{Anru}\binits{A.}}
(\byear{2013}).
\btitle{Sharp RIP bound for sparse signal and low-rank matrix recovery}.
\bjournal{Appl. Comput. Harmon. Anal.}
\bvolume{35}
\bpages{74--93}.
\bid{doi={10.1016/j.acha.2012.07.010}, issn={1063-5203}, mr={3053747}}
\end{barticle}
%
\bptok{imsref}%
\endbibitem

\bibitem{CaiZhang2}
%
\begin{barticle}[mr]
\bauthor{\bsnm{Cai},~\bfnm{T.~Tony}\binits{T.~T.}} \AND
\bauthor{\bsnm{Zhang},~\bfnm{Anru}\binits{A.}}
(\byear{2013}).
\btitle{Compressed sensing and affine rank minimization under
restricted isometry}.
\bjournal{IEEE Trans. Signal Process.}
\bvolume{61}
\bpages{3279--3290}.
\bid{doi={10.1109/TSP.2013.2259164}, issn={1053-587X}, mr={3070321}}
\end{barticle}
%
\bptok{imsref}%
\endbibitem

\bibitem{CaiZhang3}
%
\begin{barticle}[mr]
\bauthor{\bsnm{Cai},~\bfnm{T.~Tony}\binits{T.~T.}} \AND
\bauthor{\bsnm{Zhang},~\bfnm{Anru}\binits{A.}}
(\byear{2014}).
\btitle{Sparse representation of a polytope and recovery in sparse
signals and low-rank matrices}.
\bjournal{IEEE Trans. Inform. Theory}
\bvolume{60}
\bpages{122--132}.
\bid{doi={10.1109/TIT.2013.2288639}, issn={0018-9448}, mr={3150915}}
\end{barticle}
%
\bptok{imsref}%
\endbibitem

\bibitem{Supplement}
%
\begin{bmisc}[author]
{\bauthor{\bsnm{Cai},~\binits{T.}} \AND
\bauthor{\bsnm{Zhang},~\binits{A.}}}
(\byear{2014}).
\bhowpublished{Supplement to ``ROP: Matrix recovery via rank-one projections.''
DOI:\doiurl{10.1214/14-AOS1267SUPP}}.
\bptok{imsref}%
\end{bmisc}
%
\endbibitem


\bibitem{Candes-Li}
%
\begin{barticle}[mr]
\bauthor{\bsnm{Cand{\`e}s},~\bfnm{Emmanuel~J.}\binits{E.~J.}},
\bauthor{\bsnm{Li},~\bfnm{Xiaodong}\binits{X.}},
\bauthor{\bsnm{Ma},~\bfnm{Yi}\binits{Y.}} \AND
\bauthor{\bsnm{Wright},~\bfnm{John}\binits{J.}}
(\byear{2011}).
\btitle{Robust principal component analysis?}
\bjournal{J. ACM}
\bvolume{58}
\bpages{Art. 11, 37}.
\bid{doi={10.1145/1970392.1970395}, issn={0004-5411}, mr={2811000}}
\end{barticle}
%
\bptok{imsref}%
\endbibitem

\bibitem{CandesOracle}
%
\begin{barticle}[mr]
\bauthor{\bsnm{Cand{\`e}s},~\bfnm{Emmanuel~J.}\binits{E.~J.}} \AND
\bauthor{\bsnm{Plan},~\bfnm{Yaniv}\binits{Y.}}
(\byear{2011}).
\btitle{Tight oracle inequalities for low-rank matrix recovery from a
minimal number of noisy random measurements}.
\bjournal{IEEE Trans. Inform. Theory}
\bvolume{57}
\bpages{2342--2359}.
\bid{doi={10.1109/TIT.2011.2111771}, issn={0018-9448}, mr={2809094}}
\end{barticle}
%
\bptok{imsref}%
\endbibitem

\bibitem{CandesRecht}
%
\begin{barticle}[mr]
\bauthor{\bsnm{Cand{\`e}s},~\bfnm{Emmanuel~J.}\binits{E.~J.}} \AND
\bauthor{\bsnm{Recht},~\bfnm{Benjamin}\binits{B.}}
(\byear{2009}).
\btitle{Exact matrix completion via convex optimization}.
\bjournal{Found. Comput. Math.}
\bvolume{9}
\bpages{717--772}.
\bid{doi={10.1007/s10208-009-9045-5}, issn={1615-3375}, mr={2565240}}
\end{barticle}
%
\bptok{imsref}%
\endbibitem

\bibitem{CandesRIPfail}
%
\begin{barticle}[mr]
\bauthor{\bsnm{Cand{\`e}s},~\bfnm{Emmanuel~J.}\binits{E.~J.}},
\bauthor{\bsnm{Strohmer},~\bfnm{Thomas}\binits{T.}} \AND
\bauthor{\bsnm{Voroninski},~\bfnm{Vladislav}\binits{V.}}
(\byear{2013}).
\btitle{Phase{L}ift: Exact and stable signal recovery from magnitude
measurements via convex programming}.
\bjournal{Comm. Pure Appl. Math.}
\bvolume{66}
\bpages{1241--1274}.
\bid{doi={10.1002/cpa.21432}, issn={0010-3640}, mr={3069958}}
\end{barticle}
%
\bptok{imsref}%
\endbibitem

\bibitem{CandesTao}
%
\begin{barticle}[mr]
\bauthor{\bsnm{Cand{\`e}s},~\bfnm{Emmanuel~J.}\binits{E.~J.}} \AND
\bauthor{\bsnm{Tao},~\bfnm{Terence}\binits{T.}}
(\byear{2010}).
\btitle{The power of convex relaxation: Near-optimal matrix completion}.
\bjournal{IEEE Trans. Inform. Theory}
\bvolume{56}
\bpages{2053--2080}.
\bid{doi={10.1109/TIT.2010.2044061}, issn={0018-9448}, mr={2723472}}
\end{barticle}
%
\bptok{imsref}%
\endbibitem

\bibitem{Chen}
%
\begin{bmisc}[auto:parserefs-M02]
\bauthor{\bsnm{Chen},~\bfnm{Y.}\binits{Y.}},
\bauthor{\bsnm{Chi},~\bfnm{Y.}\binits{Y.}} \AND
\bauthor{\bsnm{Goldsmith},~\bfnm{A.}\binits{A.}}
(\byear{2013}).
\bhowpublished{Exact and stable covariance estimation from quadratic
sampling via convex programming.
Preprint. Available at \arxivurl{arXiv:1310.0807}.}
\end{bmisc}
%
\bptok{imsref}%
\endbibitem

\bibitem{Nowak1}
%
\begin{binproceedings}[auto:parserefs-M02]
\bauthor{\bsnm{Dasarathy},~\bfnm{G.}\binits{G.}},
\bauthor{\bsnm{Shah},~\bfnm{P.}\binits{P.}},
\bauthor{\bsnm{Bhaskar},~\bfnm{B.~N.}\binits{B.~N.}} \AND
\bauthor{\bsnm{Nowak},~\bfnm{R.}\binits{R.}}
(\byear{2012}).
\btitle{Covariance sketching}.
In \bbooktitle{50th Annual Allerton Conference on
Communication, Control, and Computing}
\bpages{1026--1033}.
\end{binproceedings}
%
\bptok{imsref}%
\endbibitem

\bibitem{Nowak2}
%
\begin{bmisc}[auto:parserefs-M02]
\bauthor{\bsnm{Dasarathy},~\bfnm{G.}\binits{G.}},
\bauthor{\bsnm{Shah},~\bfnm{P.}\binits{P.}},
\bauthor{\bsnm{Bhaskar},~\bfnm{B.~N.}\binits{B.~N.}} \AND
\bauthor{\bsnm{Nowak},~\bfnm{R.}\binits{R.}}
(\byear{2013}).
\bhowpublished{Sketching sparse matrices.
Preprint. Available at \arxivurl{arXiv:1303.6544}.}
\end{bmisc}
%
\bptok{imsref}%
\endbibitem

\bibitem{Dvijotham}
%
\begin{binproceedings}[auto:parserefs-M02]
\bauthor{\bsnm{Dvijotham},~\bfnm{K.}\binits{K.}} \AND
\bauthor{\bsnm{Fazel},~\bfnm{M.}\binits{M.}}
(\byear{2010}).
\btitle{A nullspace analysis of the nuclear norm heuristic for rank
minimization}.
In \bbooktitle{2010 IEEE International Conference on Acoustics Speech and Signal Processing (ICASSP)}
\bpages{3586--3589}.
\end{binproceedings}
%
\bptok{imsref}%
\endbibitem

\bibitem{Fan}
%
\begin{barticle}[mr]
\bauthor{\bsnm{Fan},~\bfnm{Jianqing}\binits{J.}},
\bauthor{\bsnm{Fan},~\bfnm{Yingying}\binits{Y.}} \AND
\bauthor{\bsnm{Lv},~\bfnm{Jinchi}\binits{J.}}
(\byear{2008}).
\btitle{High dimensional covariance matrix estimation using a factor model}.
\bjournal{J. Econometrics}
\bvolume{147}
\bpages{186--197}.
\bid{doi={10.1016/j.jeconom.2008.09.017}, issn={0304-4076}, mr={2472991}}
\end{barticle}
%
\bptok{imsref}%
\endbibitem

\bibitem{CVX1}
%
\begin{bmisc}[auto:parserefs-M02]
\bauthor{\bsnm{Grant},~\bfnm{M.}\binits{M.}} \AND
\bauthor{\bsnm{Boyd},~\bfnm{S.}\binits{S.}}
(\byear{2012}).
\bhowpublished{CVX: Matlab software for disciplined convex
programming, version 2.0 beta.
Available at \url{http://cvxr.com/cvx}.}
\end{bmisc}
%
\bptok{imsref}%
\endbibitem

\bibitem{CVX2}
%
\begin{bincollection}[mr]
\bauthor{\bsnm{Grant},~\bfnm{Michael~C.}\binits{M.~C.}} \AND
\bauthor{\bsnm{Boyd},~\bfnm{Stephen~P.}\binits{S.~P.}}
(\byear{2008}).
\btitle{Graph implementations for nonsmooth convex programs}.
In \bbooktitle{Recent Advances in Learning and Control (a tribute to M. Vidyasagar)}
(\beditor{\bfnm{V.}\binits{V.}~\bsnm{Blondel} \bsuffix{et~al.}}, eds.).
\bseries{Lecture Notes in Control and Inform. Sci.}
\bvolume{371}
\bpages{95--110}.
\bpublisher{Springer},
\blocation{London}.
\bid{doi={10.1007/978-1-84800-155-8_7}, mr={2409077}}
\end{bincollection}
%
\bptok{imsref}%
\endbibitem

\bibitem{Gross}
%
\begin{barticle}[auto:parserefs-M02]
\bauthor{\bsnm{Gross},~\bfnm{D.}\binits{D.}},
\bauthor{\bsnm{Liu},~\bfnm{Y.~K.}\binits{Y.~K.}},
\bauthor{\bsnm{Flammia},~\bfnm{S.~T.}\binits{S.~T.}},
\bauthor{\bsnm{Becker},~\bfnm{S.}\binits{S.}} \AND
\bauthor{\bsnm{Eisert},~\bfnm{J.}\binits{J.}}
(\byear{2010}).
\btitle{Quantum state tomography via compressed sensing}.
\bjournal{Phys. Rev. Lett.}
\bvolume{105}
\bpages{150401--150404}.
\end{barticle}
%
\bptok{imsref}%
\endbibitem

\bibitem{Johnstone}
%
\begin{barticle}[mr]
\bauthor{\bsnm{Johnstone},~\bfnm{Iain~M.}\binits{I.~M.}}
(\byear{2001}).
\btitle{On the distribution of the largest eigenvalue in principal
components analysis}.
\bjournal{Ann. Statist.}
\bvolume{29}
\bpages{295--327}.
\bid{doi={10.1214/aos/1009210544}, issn={0090-5364}, mr={1863961}}
\end{barticle}
%
\bptok{imsref}%
\endbibitem

\bibitem{KLT}
%
\begin{barticle}[mr]
\bauthor{\bsnm{Koltchinskii},~\bfnm{Vladimir}\binits{V.}},
\bauthor{\bsnm{Lounici},~\bfnm{Karim}\binits{K.}} \AND
\bauthor{\bsnm{Tsybakov},~\bfnm{Alexandre~B.}\binits{A.~B.}}
(\byear{2011}).
\btitle{Nuclear-norm penalization and optimal rates for noisy low-rank
matrix completion}.
\bjournal{Ann. Statist.}
\bvolume{39}
\bpages{2302--2329}.
\bid{doi={10.1214/11-AOS894}, issn={0090-5364}, mr={2906869}}
\end{barticle}
%
\bptok{imsref}%
\endbibitem

\bibitem{Koren}
%
\begin{barticle}[auto:parserefs-M02]
\bauthor{\bsnm{Koren},~\bfnm{Y.}\binits{Y.}},
\bauthor{\bsnm{Bell},~\bfnm{R.}\binits{R.}} \AND
\bauthor{\bsnm{Volinsky},~\bfnm{C.}\binits{C.}}
(\byear{2009}).
\btitle{Matrix factorization techniques for recommender systems}.
\bjournal{Computer}
\bvolume{42}
\bpages{30--37}.
\end{barticle}
%
\bptok{imsref}%
\endbibitem

\bibitem{Laurent}
%
\begin{barticle}[mr]
\bauthor{\bsnm{Laurent},~\bfnm{B.}\binits{B.}} \AND
\bauthor{\bsnm{Massart},~\bfnm{P.}\binits{P.}}
(\byear{2000}).
\btitle{Adaptive estimation of a quadratic functional by model selection}.
\bjournal{Ann. Statist.}
\bvolume{28}
\bpages{1302--1338}.
\bid{doi={10.1214/aos/1015957395}, issn={0090-5364}, mr={1805785}}
\end{barticle}
%
\bptok{imsref}%
\endbibitem

\bibitem{Nadler}
%
\begin{barticle}[mr]
\bauthor{\bsnm{Nadler},~\bfnm{Boaz}\binits{B.}}
(\byear{2010}).
\btitle{Nonparametric detection of signals by information theoretic
criteria: Performance analysis and an improved estimator}.
\bjournal{IEEE Trans. Signal Process.}
\bvolume{58}
\bpages{2746--2756}.
\bid{doi={10.1109/TSP.2010.2042481}, issn={1053-587X}, mr={2789420}}
\end{barticle}
%
\bptok{imsref}%
\endbibitem

\bibitem{Negahban}
%
\begin{barticle}[mr]
\bauthor{\bsnm{Negahban},~\bfnm{Sahand}\binits{S.}} \AND
\bauthor{\bsnm{Wainwright},~\bfnm{Martin~J.}\binits{M.~J.}}
(\byear{2011}).
\btitle{Estimation of (near) low-rank matrices with noise and
high-dimensional scaling}.
\bjournal{Ann. Statist.}
\bvolume{39}
\bpages{1069--1097}.
\bid{doi={10.1214/10-AOS850}, issn={0090-5364}, mr={2816348}}
\end{barticle}
%
\bptok{imsref}%
\endbibitem

\bibitem{Oymak}
%
\begin{bmisc}[auto:parserefs-M02]
\bauthor{\bsnm{Oymak},~\bfnm{S.}\binits{S.}} \AND
\bauthor{\bsnm{Hassibi},~\bfnm{B.}\binits{B.}}
(\byear{2010}).
\bhowpublished{New null space results and recovery thresholds for
matrix rank minimization.
Preprint. Available at \arxivurl{arXiv:1011.6326}.}
\end{bmisc}
%
\bptok{imsref}%
\endbibitem

\bibitem{Oymak11}
binproceedings
\begin{binproceedings}[auto:parserefs-M02]
\bauthor{\bsnm{Oymak},~\bfnm{S.}\binits{S.}},
\bauthor{\bsnm{Mohan},~\bfnm{K.}\binits{K.}},
\bauthor{\bsnm{Fazel},~\bfnm{M.}\binits{M.}} \AND
\bauthor{\bsnm{Hassibi},~\bfnm{B.}\binits{B.}}
(\byear{2011}).
\btitle{A simplified approach to recovery conditions for low-rank matrices}.
In \bbooktitle{Proc. Intl. Sympo. Information Theory (ISIT)}
\bpages{2318--2322}.
\bpublisher{IEEE}, \blocation{Piscataway, NJ}.
\end{binproceedings}
%
\bptok{imsref}%
\endbibitem

\bibitem{Patterson}
%
\begin{barticle}[pbm]
\bauthor{\bsnm{Patterson},~\bfnm{Nick}\binits{N.}},
\bauthor{\bsnm{Price},~\bfnm{Alkes~L.}\binits{A.~L.}} \AND
\bauthor{\bsnm{Reich},~\bfnm{David}\binits{D.}}
(\byear{2006}).
\btitle{Population structure and eigenanalysis}.
\bjournal{PLoS Genet.}
\bvolume{2}
\bpages{e190}.
\bid{doi={10.1371/journal.pgen.0020190}, issn={1553-7404},
pii={06-PLGE-RA-0101R3}, pmcid={1713260}, pmid={17194218}}
\end{barticle}
%
\bptok{imsref}%
\endbibitem

\bibitem{Price}
%
\begin{barticle}[pbm]
\bauthor{\bsnm{Price},~\bfnm{Alkes~L.}\binits{A.~L.}},
\bauthor{\bsnm{Patterson},~\bfnm{Nick~J.}\binits{N.~J.}},
\bauthor{\bsnm{Plenge},~\bfnm{Robert~M.}\binits{R.~M.}},
\bauthor{\bsnm{Weinblatt},~\bfnm{Michael~E.}\binits{M.~E.}},
\bauthor{\bsnm{Shadick},~\bfnm{Nancy~A.}\binits{N.~A.}} \AND
\bauthor{\bsnm{Reich},~\bfnm{David}\binits{D.}}
(\byear{2006}).
\btitle{Principal components analysis corrects for stratification in
genome-wide association studies}.
\bjournal{Nat. Genet.}
\bvolume{38}
\bpages{904--909}.
\bid{doi={10.1038/ng1847}, issn={1061-4036}, pii={ng1847}, pmid={16862161}}
\end{barticle}
%
\bptok{imsref}%
\endbibitem

\bibitem{RechtImproved}
%
\begin{barticle}[mr]
\bauthor{\bsnm{Recht},~\bfnm{Benjamin}\binits{B.}}
(\byear{2011}).
\btitle{A simpler approach to matrix completion}.
\bjournal{J. Mach. Learn. Res.}
\bvolume{12}
\bpages{3413--3430}.
\bid{issn={1532-4435}, mr={2877360}}
\end{barticle}
%
\bptok{imsref}%
\endbibitem

\bibitem{RechtMatrix}
%
\begin{barticle}[mr]
\bauthor{\bsnm{Recht},~\bfnm{Benjamin}\binits{B.}},
\bauthor{\bsnm{Fazel},~\bfnm{Maryam}\binits{M.}} \AND
\bauthor{\bsnm{Parrilo},~\bfnm{Pablo~A.}\binits{P.~A.}}
(\byear{2010}).
\btitle{Guaranteed minimum-rank solutions of linear matrix equations
via nuclear norm minimization}.
\bjournal{SIAM Rev.}
\bvolume{52}
\bpages{471--501}.
\bid{doi={10.1137/070697835}, issn={0036-1445}, mr={2680543}}
\end{barticle}
%
\bptok{imsref}%
\endbibitem

\bibitem{Rohde}
%
\begin{barticle}[mr]
\bauthor{\bsnm{Rohde},~\bfnm{Angelika}\binits{A.}} \AND
\bauthor{\bsnm{Tsybakov},~\bfnm{Alexandre~B.}\binits{A.~B.}}
(\byear{2011}).
\btitle{Estimation of high-dimensional low-rank matrices}.
\bjournal{Ann. Statist.}
\bvolume{39}
\bpages{887--930}.
\bid{doi={10.1214/10-AOS860}, issn={0090-5364}, mr={2816342}}
\end{barticle}
%
\bptok{imsref}%
\endbibitem

\bibitem{Trosset}
%
\begin{barticle}[mr]
\bauthor{\bsnm{Trosset},~\bfnm{Michael~W.}\binits{M.~W.}}
(\byear{2000}).
\btitle{Distance matrix completion by numerical optimization}.
\bjournal{Comput. Optim. Appl.}
\bvolume{17}
\bpages{11--22}.
\bid{doi={10.1023/A:1008722907820}, issn={0926-6003}, mr={1791595}}
\end{barticle}
%
\bptok{imsref}%
\endbibitem

\bibitem{Vershynin}
%
\begin{barticle}[mr]
\bauthor{\bsnm{Vershynin},~\bfnm{Roman}\binits{R.}}
(\byear{2011}).
\btitle{Spectral norm of products of random and deterministic matrices}.
\bjournal{Probab. Theory Related Fields}
\bvolume{150}
\bpages{471--509}.
\bid{doi={10.1007/s00440-010-0281-z}, issn={0178-8051}, mr={2824864}}
\end{barticle}
%
\bptok{imsref}%
\endbibitem

\bibitem{Wakinimage}
%
\begin{binproceedings}[auto:parserefs-M02]
\bauthor{\bsnm{Wakin},~\bfnm{M.}\binits{M.}},
\bauthor{\bsnm{Laska},~\bfnm{J.}\binits{J.}},
\bauthor{\bsnm{Duarte},~\bfnm{M.}\binits{M.}},
\bauthor{\bsnm{Baron},~\bfnm{D.}\binits{D.}},
\bauthor{\bsnm{Sarvotham},~\bfnm{S.}\binits{S.}},
\bauthor{\bsnm{Takhar},~\bfnm{D.}\binits{D.}},
\bauthor{\bsnm{Kelly},~\bfnm{K.}\binits{K.}} \AND
\bauthor{\bsnm{Baraniuk},~\bfnm{R.}\binits{R.}}
(\byear{2006}).
\btitle{An architecture for compressive imaging}.
In \bbooktitle{Proceedings of the International Conference on Image
Processing (ICIP 2006)}
\bpages{1273--1276}.
\end{binproceedings}
%
\bptok{imsref}%
\endbibitem

\bibitem{WangLi}
%
\begin{barticle}[auto:parserefs-M02]
\bauthor{\bsnm{Wang},~\bfnm{H.}\binits{H.}} \AND
\bauthor{\bsnm{Li},~\bfnm{S.}\binits{S.}}
(\byear{2013}).
\btitle{The bounds of restricted isometry constants for low rank matrices recovery}.
\bjournal{Sci. China Ser. A}
\bvolume{56}
\bpages{1117--1127}.
\end{barticle}
%
\bptok{imsref}%
\endbibitem

\bibitem{Wang}
%
\begin{barticle}[mr]
\bauthor{\bsnm{Wang},~\bfnm{Yazhen}\binits{Y.}}
(\byear{2013}).
\btitle{Asymptotic equivalence of quantum state tomography and noisy
matrix completion}.
\bjournal{Ann. Statist.}
\bvolume{41}
\bpages{2462--2504}.
\bid{doi={10.1214/13-AOS1156}, issn={0090-5364}, mr={3127872}}
\end{barticle}
%
\bptok{imsref}%
\endbibitem

\bibitem{Wax}
%
\begin{barticle}[mr]
\bauthor{\bsnm{Wax},~\bfnm{Mati}\binits{M.}} \AND
\bauthor{\bsnm{Kailath},~\bfnm{Thomas}\binits{T.}}
(\byear{1985}).
\btitle{Detection of signals by information theoretic criteria}.
\bjournal{IEEE Trans. Acoust. Speech Signal Process.}
\bvolume{33}
\bpages{387--392}.
\bid{doi={10.1109/TASSP.1985.1164557}, issn={0096-3518}, mr={0788604}}
\end{barticle}
%
\bptok{imsref}%
\endbibitem

\end{thebibliography}
\end{document}